\documentclass[leqno,12pt]{article}
\usepackage{amsmath,amssymb,rotating,a4wide,color}
\usepackage{wasysym}
\usepackage[all,cmtip]{xy}

\usepackage{verbatim}

\newdir^{ (}{{}*!/-3pt/\dir^{(}}
\newdir_{ (}{{}*!/-3pt/\dir_{(}}

\entrymodifiers={+!!<0pt,\fontdimen22\textfont2>}

\def\opp{{\rm opp}}
\def\tTp{{\mathstrut\smash{\widetilde{T_\Fp}}}}

\def\oK{{\overline{K}}}
\def\oL{{\overline{L}}}
\def\UX{{\underline{X}}}
\def\smalldash{{\raisebox{1pt}{$\scriptscriptstyle-$}}}

\def\bigtimes{\mathop{\raise-2pt\hbox{\huge$\times$}}}
\def\Ux{{\underline{x}}}
\def\GammaUx{\Gamma_{\kern-1pt\Ux}}

\newbox\circbulletbox
\setbox\circbulletbox=\hbox{\,{\large$\circledcirc$}\kern-8.7pt\raise .6pt\hbox{$\bullet$}\,}

\let\le\leqslant
\let\ge\geqslant
\let\leq\leqslant
\let\geq\geqslant

\def\mycirc{{\kern1pt\circ\kern2pt}}
\def\charact{\mathop{\rm char}\nolimits}

\def\Image{\mathop{\rm Im}\nolimits}

\def\Aut{\mathop{\rm Aut}\nolimits}
\def\Frob{\mathop{\rm Frob}\nolimits}
\def\Hom{\mathop{\rm Hom}\nolimits}

\def\Gal{\mathop{\rm Gal}\nolimits}
\def\End{\mathop{\rm End}\nolimits}

\def\Spec{\mathop{\rm Spec}\nolimits}

\def\deg{\mathop{\rm deg}\nolimits}

\def\Ker{\mathop{\rm Ker}\nolimits}

\def\Mat{\mathop{\rm Mat}\nolimits}

\def\Cent{\mathop{\rm Cent}\nolimits}

\def\Quot{\mathop{\rm Quot}\nolimits}

\def\trace{\mathop{\rm tr}\nolimits}

\def\ord{\mathop{\rm ord}\nolimits}
\def\rank{\mathop{\rm rank}\nolimits}

\def\GL{\mathop{\rm GL}\nolimits}
\def\SL{\mathop{\rm SL}\nolimits}

\def\ad{{\rm ad}}

\def\trad{{\rm trad}}

\def\sep{{\rm sep}}

\let\phi\varphi
\let\theta\vartheta
\let\epsilon\varepsilon
\let\setminus\smallsetminus

\let\emptyset\varnothing

\newtheorem{Thm}{Theorem}[section]
\newtheorem{Prop}[Thm]{Proposition}
\newtheorem{Lem}[Thm]{Lemma}

\newtheorem{Def}[Thm]{Definition}

\newtheorem{Rem}[Thm]{Remark}

\newtheorem{Var}[Thm]{Variation}

\numberwithin{Thm}{section}
\def\UseTheoremCounterForNextEquation{\setcounter{equation}{\value{Thm}}\addtocounter{Thm}{1}}

\def\qed{{\hskip0pt\unskip\unskip\nobreak\hfil\penalty50
          \hskip1em\hbox{}\nobreak\hfil
           {$\square$}
          \parfillskip=0pt\finalhyphendemerits=0
          \par}\medskip}

\newenvironment{Proof}
               {\noindent{\bf Proof.}\ }
               {\qed}

               {\noindent{\bf Proof of #1.}\ }
               {\qed}

\newcommand{\BF}{{\mathbb{F}}}
\newcommand{\BG}{{\mathbb{G}}}

\newcommand{\BZ}{{\mathbb{Z}}}

\newcommand{\Fa}{{\mathfrak{a}}}

\newcommand{\Fm}{{\mathfrak{m}}}
\newcommand{\Fn}{{\mathfrak{n}}}

\newcommand{\Fp}{{\mathfrak{p}}}
\newcommand{\Fq}{{\mathfrak{q}}}

\newcommand{\FP}{{\mathfrak{P}}}

\newcommand{\CL}{{\cal L}}

\newcommand{\CO}{{\cal O}}

\newcommand{\CR}{{\cal R}}
\newcommand{\CS}{{\cal S}}

\newbox\mybox
\def\arrover#1{\mathrel{
       \setbox\mybox=\hbox spread 1.4em
              {\hfil$\scriptstyle#1$\hfil}
       \vbox{\offinterlineskip\copy\mybox
             \hbox to\wd\mybox{\rightarrowfill}}}}

\def\larrover#1{\mathrel{
       \setbox\mybox=\hbox spread 1.4em
              {\hfil$\scriptstyle#1\vphantom{g}$\hfil}
       \vbox{\offinterlineskip\copy\mybox
             \hbox to\wd\mybox{\leftarrowfill}}}}

\def\ontoover#1{\mathrel{
       \setbox\mybox=\hbox spread 1.4em
              {\hfil$\scriptstyle#1\vphantom{g}$\hfil}
       \vbox{\offinterlineskip\copy\mybox
             \hbox to\wd\mybox{\rightarrowfill\hskip-2.8mm
                               $\rightarrow$}}}}
\def\leftontoover#1{\mathrel{
       \setbox\mybox=\hbox spread 1.4em
              {\hfil$\scriptstyle#1\vphantom{g}$\hfil}
       \vbox{\offinterlineskip\copy\mybox
             \hbox to\wd\mybox{$\leftarrow$\hskip-2.8mm
                               \leftarrowfill}}}}
\let\longto\longrightarrow
\let\into\hookrightarrow
\let\onto\twoheadrightarrow

\begin{document}

\title{Finding Endomorphisms of Drinfeld modules}

\author{\begin{tabular}{cc}
Nikolas Kuhn
& Richard Pink\\[12pt]
  \small Dept. of Mathematics
& \small Dept. of Mathematics \\
  \small ETH Z\"urich 
& \small ETH Z\"urich\\
% \small 8092 Z\"urich\\
  \small Switzerland
& \small Switzerland\\
\small \tt nikolaskuhn@gmx.de
& \small \tt pink@math.ethz.ch
\end{tabular}}

%\author{Richard Pink\\[12pt]
%\small Department of Mathematics \\[-3pt]
%\small ETH Z\"urich\\[-3pt]
%\small 8092 Z\"urich\\[-3pt]
%\small Switzerland \\[-3pt]
%\small pink@math.ethz.ch\\[12pt]}

\date{August 9, 2016}

\maketitle

\begin{abstract}
We give an effective algorithm to determine the endomorphism ring of a Drinfeld module, both over its field of definition and over a separable or algebraic closure thereof. Using previous results we deduce an effective description of the image of the adelic Galois representation associated to the Drinfeld module, up to commensurability. 
We also give an effective algorithm to decide whether two Drinfeld modules are isogenous, again both over their field of definition and over a separable or algebraic closure thereof.
\end{abstract}

{\renewcommand{\thefootnote}{}
\footnotetext{
% 11G09 Drinfeld modules; higher-dimensional motives, etc.
% 11R58 Arithmetic theory of algebraic function fields
% 11Y99 Computational number theory 
% Keywords:
% Drinfeld modules, endomorphisms, Galois representations, effective algorithms
}
}

%\newpage
% This reduces the line spacing in the table of contents:
\renewcommand{\baselinestretch}{0.75}\normalsize
\tableofcontents
\renewcommand{\baselinestretch}{1.0}\normalsize
\newpage

%%%%%%%%%%%%%%%%%%%%%%%%%%%%%%%%%%%%%%%%%%%%%%%%%%%%%%%%%%%%%%%%%%%%%%%%%%%%
%%%%%%%%%%%%%%%%%%%%%%%%%%%%%%%%%%%%%%%%%%%%%%%%%%%%%%%%%%%%%%%%%%%%%%%%%%%%

\section{Introduction}
\label{Intro}

Given a Drinfeld $A$-module $\phi\colon A\to K[\tau]$ over a field~$K$, can one effectively determine its endomorphism ring $\End_K(\phi)$?

\medskip
Before answering this question, we must make it more precise. By definition $\End_K(\phi)$ is the subset of elements of $K[\tau]$ which commute with $\phi_a$ for all $a\in A$. Thus one can write down individual endomorphisms, but what does it mean to know their totality? We think it means three things: 
Firstly, since $\End_K(\phi)$ is a finitely generated projective $A$-module, one should have a finite set of generators. Secondly, one should know all $A$-linear relations between them, in other words, one should have a finite presentation of $\End_K(\phi)$ as an $A$-module. Thirdly, one should be able to express any given endomorphism as an $A$-linear combination of the generators. Applying this to the product of any two generators, this yields an explicit description of the ring structure of $\End_K(\phi)$, using which many questions about $\End_K(\phi)$ as an $A$-algebra reduce to finite calculations over~$A$. 

To answer the question we must also clarify which algebraic calculations we assume that one can already perform within $A$ and within~$K$. 
By definition, the coefficient ring $A$ underlying a Drinfeld module is a finitely generated normal integral domain of transcendence degree $1$ over a finite prime field $\BF_p$ of order~$p$. So we assume that $A$ is given by explicit finite sets of generators and relations. 
We also assume that $K$ is the fraction field of an integral domain that is given by explicit finite sets of generators and relations over~$\BF_p$. For calculations within $A$ and $K$ we then have all the standard procedures from algorithmic commutative algebra at our disposal. 
Of course, one cannot effectively construct, or calculate within, a separable or algebraic closure $K^\sep\subset\overline{K}$ of~$K$. 
But one can calculate in any finite extension of $K$ and enlarge that extension whenever necessary. 
%Throughout all finite separable extensions of $K$ are tacitly assumed to be contained in~$K^\sep$.

The assumption that $K$ is finitely generated over~$\BF_p$, however, introduces a new problem. Namely, while 
%$\End_{\oK}(\phi)$ is still a finitely generated $A$-module of bounded rank, and hence
there exists a finite field extension $K'$ of $K$ 
%of bounded degree 
with $\End_{\oK}(\phi)=\End_{K'}(\phi)$, there is no a priori choice for it.
%\Bigskip (NOTE: One can first choose a prime $\Fp\not=\Fp_0$ of~$A$ and explicitly find a finite separable extension $K'$ of $K$ such that all $\Fp$-torsion points are defined over~$K'$. Then by looking at the action of endomorphisms on the $\Fp$-torsion one can probably say something about their fields of definition, at least for endomorphisms of degree prime to~$\Fp$.) \Bigskip
To determine $\End_{\oK}(\phi)$ we must therefore also specify such an extension~$K'$.

\medskip
With these provisos we can now say that $\End_K(\phi)$ and  $\End_{\oK}(\phi)$ can be effectively determined: see Theorems \ref{AllEndosK} and \ref{AllEndosKSep}.

\medskip
Our algorithm for this has essentially two parts. One process goes through all integers $d\ge0$ and finds all endomorphisms of degree $d$ by solving finitely many polynomial equations. Eventually it will find a finite set of generators, but knowing when that happens requires other information. It is not hard to see that it suffices to know the rank of the endomorphism ring over~$A$.
So in addition to the first process, we start another process in parallel that tries to prove that the right number of $A$-linearly independent endomorphisms has already been found. When that succeeds, it kills the first process and stops with the correct answer. 

The second process uses the Galois representation on the $\Fp$-adic Tate module of $\phi$ for a suitable prime $\Fp$ of~$A$. By the Tate conjecture for Drinfeld modules, proved by Taguchi \cite{TaguchiFinite}, \cite{TaguchiInfinite}, \cite{TaguchiTate} and Tamagawa \cite{Tamagawa0}, this representation determines the endomorphism ring to a large extent; in particular, it determines its rank over~$A$. 
Though the Galois representation can be studied only indirectly, the associated characteristic polynomials of Frobeniuses can be computed effectively.
The second process limits the possible endomorphisms by searching for characteristic polynomials that are sufficiently independent of each other in some sense.

\medskip
For Drinfeld modules of generic characteristic, the endomorphism ring is always commutative, and the program outlined above suffices to determine it effectively. In special characteristic the endomorphism ring can be non-commutative, and we must wrestle with additional technical difficulties. The problem is that there may exist more endomorphisms when $A$ is replaced by a smaller admissible coefficient ring. This puts additional constraints on the Galois representation. 

\medskip
In fact, by results of the second author and others \cite{PinkRIMS}, \cite{PinkDrinSpec2}, \cite{PinkTraulsen1}, \cite{PinkRuetsche2}, \cite{DevicPink}, the knowledge of the endomorphism rings of certain Drinfeld modules obtained from $\phi$ by varying the ring of coefficients~$A$ determines the image of the Galois representation up to commensurability, even for the whole adelic Galois representation associated to~$\phi$. We therefore set up things to compute all this information as well and are thereby able to effectively determine the image of Galois up to commensurability: see Theorem \ref{GaloisUpToCommens}.

\medskip
In special characteristic we treat the isotrivial case, where $\phi$ is isomorphic over $\oK$ to a Drinfeld module defined over a finite field, separately. In the non-isotrivial case we first use the method sketched above to find a maximal commutative subring $A'$ of $\End_{\oK}(\phi)$. For technical reasons we replace $\phi$ by an isogenous Drinfeld module, after which $A'$ is an admissible coefficient ring and $\phi$ extends to a Drinfeld $A'$-module with $\End_{\oK}(\phi')=A'$. In a second step we then find the unique smallest admissible coefficient ring $B\subset A'$ such that the center of $\End_{K^\sep}(\phi'|B)$ is~$B$, whose existence is guaranteed by Pink \cite[Thm.\;1.2]{PinkDrinSpec2}.

Our algorithm for this again has two parts. One process computes the traces of Frobeniuses in the adjoint representation, whose values generate the fraction field of $B$ by Pink \cite[Thm.\;1.3]{PinkDrinSpec2}. It thus constructs an increasing sequence of subrings $B_k$ of~$A'$ with $B_k=B$ for all sufficiently large~$k$, but again it does not know when that occurs. In addition to the first process, we therefore run another process in parallel that tries to prove that $B$ has been reached. This process is started as soon as the first $B_k$ is infinite, and it simply goes through all integers $d\ge0$ and finds all endomorphisms of degree $d$ of $\phi'|B_k$ over~$\oK$. This process knows whether $B_k=B$ has been reached by computing the ranks of $B_k$ and of the submodule of $\End_{\oK}(\phi'|B_k)$ that is generated by the endomorphisms already found. When that occurs, it kills the first process and stops with the correct answer. 
Using the knowledge of $B$ one can then find the endomorphism ring of the original Drinfeld module~$\phi$.

\medskip
In all this, we do not care about computational efficiency; we only try to keep the code short and well-organized. 
An actual implementation should probably retain intermediate information and reuse it later. Also, simply searching for all endomorphisms of small degree seems a brute force approach. In Section \ref{Var} we discuss some ideas which might speed up the search by finding a priori candidates for the generators of the endomorphism ring.

\medskip
A natural related question, kindly raised by Peter Jossen, is whether one can effectively decide whether two given Drinfeld $A$-modules $\phi$ and $\psi$ over $K$ are isogenous over~$K$, respectively over~$\oK$. Using the same methods as for endomorphisms, we answer this question affirmatively and show that $\Hom_K(\phi,\psi)$ and $\Hom_{\oK}(\phi,\psi)$ can be effectively determined: see Section \ref{Compare}.

\medskip
With effectiveness established, one may ask whether there exist any kinds of a priori bounds on the endomorphism ring, say on its rank and its discriminant over~$A$, for instance in terms of the height of the Drinfeld module. (For abelian varieties over number fields such bounds are due to Masser and W\"ustholz \cite{MasserWuestholz1994}.) This article has nothing to contribute to this question, but it may be an interesting one for someone to pursue in the future.

Also, we have not tried to actually implement the proposed algorithms and can therefore not show any nice examples.

\medskip
The article grew out of the master thesis of the first author \cite{KuhnMasterThesis}.

%Section \ref{EndGal}: Endomorphisms and image of Galois
%Section \ref{Comp}: Computer algebra prerequisites
%Section \ref{AlgBits}: Bits of algorithms
%Section \ref{SearchEndos}: Searching for endomorphisms
%Section \ref{MainAlg}: Main algorithms
%Section \ref{Var}: Variation
%Section \ref{Compare}: Comparing two Drinfeld modules

%%%%%%%%%%%%%%%%%%%%%%%%%%%%%%%%%%%%%%%%%%%%%%%%%%%%%%%%%%%%%%%%%%%%%%%%%%%
%\newpage
\section{Endomorphisms and image of Galois}
\label{EndGal}

In this section we review known facts about endomorphisms and Galois representations associated to Drinfeld modules and deduce some consequences. For the general theory of Drinfeld modules see Drinfeld \cite{Drinfeld1}, Deligne and Husem\"oller \cite{DelHus}, Hayes \cite{HayEx}, or Goss \cite{GossBS}.

\medskip
{\bf Basics:} Let $\BF_p$ denote the finite field of prime order~$p$. Let $F$ be a finitely generated field of transcendence degree $1$ over~$\BF_p$, and let $A$ be the subring of elements of $F$ which are regular outside a fixed place $\infty$ of~$F$. We call such $A$ an admissible coefficient ring.

Let $K$ be another finitely generated field over~$\BF_p$ with separable, respectively algebraic closures $K^\sep\subset\oK$.
% Let $K[\tau]$ denote the non-commutative polynomial ring in one variable over~$K$, where $\tau$ satisfies the commutation relation $\tau u=u^p\tau$ for all $u\in K$.
Write $\End(\BG_{a,K}) = K[\tau]$ with $\tau(x) = x^p$. 
Consider a Drinfeld $A$-module ${\phi\colon A \to K[\tau]}$, $a\mapsto\phi_a$ of rank~$r$ with characteristic ideal $\Fp_0\subset A$.
%Its  characteristic ideal $\Fp_0$ is the kernel of the homomorphism $\phi\colon A\to K$ determined by the lowest coefficient of~$\phi$. 
Recall that $\phi$ has generic characteristic if $\Fp_0=(0)$ and special characteristic otherwise. We call $\phi$ isotrivial if, over some finite extension of~$K$,  it is isomorphic to a Drinfeld $A$-module defined over a finite field; this can happen only in special characteristic.

\medskip
{\bf Endomorphisms:} By definition $\End_K(\phi)$ is the centralizer of $\phi(A)$ in $K[\tau]$. This is a finitely generated projective $A$-module, and $\End_K^\circ(\phi) := \End_K(\phi)\otimes_AF$ is a division algebra over~$F$ of finite dimension dividing~$r^2$. There exists a subfield $K'\subset K^\sep$ finite over~$K$ such that $\End_\oK(\phi) = \End_{K'}(\phi)$.
% There exists a finite separable extension $K'$ of $K$ such that for every overfield $L$ of $K'$ we have $\End_L(\phi) = \End_{K'}(\phi)$.
In generic characteristic the endomorphism ring is always commutative. 

\medskip
{\bf Good reduction:} Choose a normal integral domain $R\subset K$ which is finitely generated over~$\BF_p$ with $\Quot(R)=K$, such that $\phi$ extends to a Drinfeld $A$-module over $\Spec R$. For any maximal ideal $\Fm\subset R$ let $\phi_\Fm$ denote the resulting Drinfeld $A$-module over the finite residue field $k_\Fm:=R/\Fm$. It is known that any endomorphism of $\phi$ over $K$ already has coefficients in~$R$; so reduction modulo $\Fm$ induces a natural homomorphism of $A$-algebras
\UseTheoremCounterForNextEquation
\begin{equation}\label{EndRed}
\End_K(\phi) \longto \End_{k_\Fm}(\phi_\Fm).
\end{equation}
Moreover, the degree in $\tau$ of an endomorphism is preserved under reduction; hence the homomorphism is injective.

\medskip
{\bf Frobenius:} The element $\Frob_\Fm := \tau^{[k_\Fm/\BF_p]}$ lies in the center of $k_\Fm[\tau]$ and therefore in $\End_{k_\Fm}(\phi_\Fm)$. 
In fact, the center of $\End_{k_\Fm}^\circ(\phi_\Fm)$ is the field extension $F(\Frob_\Fm)$ of $F$ that is generated by $\Frob_\Fm$. 
Moreover, let $d_\Fm$ denote the dimension of $F(\Frob_\Fm)$ over~$F$, and let $e_\Fm^2$ be the dimension of $\End_{k_\Fm}^\circ(\phi_\Fm)$ over $F(\Frob_\Fm)$; then we have $d_\Fm e_\Fm=r$.

Let $\min_\Fm(X)$ denote the minimal polynomial of $\Frob_\Fm$ over~$F$; by construction it is irreducible and monic of degree~$d_\Fm$. Since $\Frob_\Fm$ lies in an $A$-algebra of finite rank, this polynomial actually has coefficients in~$A$. 
Define $\charact_\Fm(X) := \min_\Fm(X)^{e_\Fm}$, which is a monic polynomial in $A[X]$ of degree~$r$, called the characteristic polynomial of $\Frob_\Fm$. 

%Let $k_\infty$ denote the finite residue field of the place $\infty$ of~$F$. Let $\oF$ be an algebraic closure of $F$, and let $\oinfty$ be a place of $\oF$ lying over~$\infty$. Let $\ord_\oinfty$ denote the associated valuation on $\oF$ whose restriction to $F$ is normalized. Then any zero $\alpha\in\oF$ of the characteristic polynomial $\charact_\Fm(X)$ is integral over~$A$ and satisfies
%% Drinfeld \cite[Prop.\;2.1]{Drinfeld2}
%\UseTheoremCounterForNextEquation
%\begin{equation}\label{EigVal}
%\ord_{\overline{\infty}}(\alpha)
%\ =\ -\frac{1}{r}\cdot\frac{[k_\Fm/\BF_p]}{[k_\infty/\BF_p]}.
%\end{equation}

\medskip
{\bf Tate modules:} For any maximal ideal $\Fp\not=\Fp_0$ of $A$ the $\Fp$-adic Tate module $T_\Fp(\phi)$ is a free module of rank $r$ over the completion~$A_\Fp$.
% which is naturally associated to $\phi$ and $\Fp$ and the choice of $K^\sep$.
It is naturally endowed with an action of $\End_K(\phi)$ and a continuous action of the Galois group $\Gal(K^\sep/K)$. These actions commute with each other, and each helps in understanding the other.

Let $R^\sep$ denote the integral closure of $R$ in~$K^\sep$. For any maximal ideal $\Fm\subset R$ choose a maximal ideal $\Fm^\sep\subset R^\sep$ which contains~$\Fm$. Then its residue field $k_\Fm^\sep:=R^\sep/\Fm^\sep$ is a separable closure of~$k_\Fm$. For any maximal ideal $\Fp$ of $A$ different from the characteristic ideal of $\phi_\Fm$ this choice induces a natural isomorphism $T_\Fp(\phi) \cong T_\Fp(\phi_\Fm)$. This isomorphism is compatible with the action of endomorphisms via the reduction homomorphism (\ref{EndRed}). It is also compatible with the action of the decomposition group at $\Fm^\sep$; namely, the inertia group acts trivially on $T_\Fp(\phi)$, and the isomorphism is equivariant under the action of the Frobenius at~$\Fm$. Moreover, the characteristic polynomial of this Frobenius in its action on the Tate module is precisely the characteristic polynomial $\charact_\Fm(X)$ defined above.

\medskip
{\bf Adelic Galois representation:} The product $T_\ad(\phi) = \prod_{\Fp\not=\Fp_0}T_\Fp(\phi)$ is a free module of rank $r$ over $A_\ad = \prod_{\Fp\not=\Fp_0} A_\Fp$, called the prime-to-$\Fp_0$ adelic Tate module of~$\phi$. It again carries natural commuting actions of $\End_K(\phi)$ and of $\Gal(K^\sep/K)$. The latter corresponds to a continuous homomorphism 
\UseTheoremCounterForNextEquation
\begin{equation}\label{GamadDef}
\rho_\ad\colon \Gal(K^\sep/K) \ \to\ \Aut_{A_\ad}(T_\ad(\phi)) \ \cong\ \GL_r(A_\ad).
\end{equation}
The image of $\rho_\ad$ is determined up to commensurability by endomorphisms, as described below.

\medskip
{\bf Isogenies:} A non-zero homomorphism between two Drinfeld $A$-modules is called an isogeny. With an isogeny, we can often reduce ourselves to Drinfeld $A'$-modules of smaller rank for a larger ring~$A'$, using the following fact: 

\begin{Prop}\label{Isog}
(Hayes \cite[Prop.\;3.2]{HayEx}, Devic-Pink \cite[Prop.\;4.3]{DevicPink})
Let $A^\smalldash$ be a commutative $A$-subalgebra of $\End_K(\phi)$. Then its normalization $A'$ is an admissible coefficient ring, and there exist a Drinfeld $A'$-module $\phi'\colon A'\to K[\tau]$ and an isogeny $h\colon\phi\to\phi'|A$ over~$K$. Moreover, we have $\rank(\phi) = \rank_A(A')\cdot\rank(\phi')$.
\end{Prop}

For the remainder of the present section we fix a maximal commutative subring $A^\smalldash$ of $\End_{K^\sep}(\phi)$ and a subfield $K'\subset K^\sep$ which is finite over~$K$ such that $A^\smalldash\subset\End_{K'}(\phi)$. Set $F':=\Quot(A^\smalldash)$ and let $A'\subset F'$ be the normalization of~$A^\smalldash$. Using Proposition \ref{Isog} over~$K'$ we choose a Drinfeld $A'$-module $\phi'\colon A'\to K'[\tau]$ and an isogeny $h\colon\phi\to\phi'|A$ over~$K'$. Then $\End_{K^\sep}(\phi')=A'$. Moreover $\phi$ is of special characteristic, respectively isotrivial, if and only if $\phi'$ is so.
Set $r':=\rank(\phi')$ and consider the adelic Galois representation associated to~$\phi'$:
\UseTheoremCounterForNextEquation
\begin{equation}\label{Gamad'Def}
\rho'_\ad\colon \Gal(K^\sep/K') \ \to\ \Aut_{A'_\ad}(T_\ad(\phi')) \ \cong\ \GL_{r'}(A'_\ad).
\end{equation}

\medskip
{\bf Generic characteristic:} Here the image of Galois is described by:

\begin{Thm}\label{EndGalGen}
(Pink-R\"utsche \cite{PinkRuetsche2})
If $\phi$ has generic characteristic, the image of $\rho'_\ad$ is an open subgroup of 
$\GL_{r'}(A'_\ad)$, and the image of $\rho_\ad$ is commensurable with the subgroup
$$\Cent_{\GL_r(A_\ad)}(\End_{K^\sep}(\phi)).$$
\end{Thm}

\medskip
{\bf Special characteristic:} Here the endomorphism ring may be non-commutative; moreover, there may exist an admissible coefficient ring $B\subsetneqq A$ with $\End_{K^\sep}(\phi)\subsetneqq\End_{K^\sep}(\phi|B)$, which puts additional constraints on the image of Galois. If $\phi$ is isotrivial, the image of $\rho_\ad$ is commensurable with the pro-cyclic subgroup generated by the image of Frobenius. Otherwise:

\begin{Thm}\label{EndGalSpec1}
(Pink \cite[Thm.\;1.2]{PinkDrinSpec2})
If $\phi$ is non-isotrivial of special characteristic, there exists a unique admissible coefficient ring $B\subset A'$ with the properties:
\begin{enumerate}
\item[(a)] The center of $\End_{K^\sep}(\phi'|B)$ is~$B$.
\item[(b)] For every admissible coefficient ring $B'\subset A'$ we have $\End_{K^\sep}(\phi'|B') \subset \End_{K^\sep}(\phi'|B)$.
\end{enumerate}
\end{Thm}

In almost all cases this subring $B$ can be characterized independently using traces of Frobenius. In fact $B$ is determined by the subfield $E :=\Quot(B)$ of $F':=\Quot(A')$, because $B=A'\cap E$.
Choose a normal integral domain $R'\subset K'$ which is finitely generated over~$\BF_p$ with $\Quot(R')=K'$, such that $\phi'$ extends to a Drinfeld $A'$-module over $\Spec R'$. For any maximal ideal $\Fm'\subset R'$ let $\phi'_{\Fm'}$ denote the resulting Drinfeld $A'$-module over the finite residue field $k_{\Fm'}:=R'/\Fm'$. 
Write the characteristic polynomial of $\Frob_{\Fm'}$ associated to $\phi'$ in the form $\sum_{i=0}^{r'}a_iX^i$ with $a_i\in F'$, or in the form $\prod_{i=1}^{r'}(X-\alpha_i)$ over an algebraic closure of~$F'$, and set 
\UseTheoremCounterForNextEquation
\begin{equation}\label{tmDef}
t_{\Fm'}\ :=\ \frac{a_1a_{r'-1}}{a_0} 
\ =\ \smash{\sum_{i=1}^{\;r'}\sum_{j=1}^{\;r'}\frac{\alpha_i}{\alpha_j}}
\ \in\ F'.
\end{equation}
Thus $t_{\Fm'}$ is the trace of $\Frob_{\Fm'}$ in the adjoint representation on $\End_{A'_{\Fp'}}(T_{\Fp'}(\phi'))$
% $T_{\Fp'}(\phi')\otimes_{A'_{\Fp'}}T_{\Fp'}(\phi')^\vee$
for any maximal ideal $\Fp'$ different from the characteristic of~$\phi'_{\Fm'}$. Let $E^\trad\subset F'$ be the subfield generated by the elements $t_{\Fm'}$ for all~$\Fm'$.

\begin{Thm}\label{EndGalSpec2}
(Pink \cite[Thm.\;1.3]{PinkDrinSpec2})
In the situation of Theorem \ref{EndGalSpec1}, we have either
\begin{enumerate}
\item[(c)] $E^\trad=E$, or
\item[(c')] 
$p=\rank(\phi')=2$ and $E^\trad=\{e^2\mid e\in E\}$.
\end{enumerate}
\end{Thm}

To describe the image of Galois let $r''$ be the rank and $\Fq_0\subset B$ the characteristic ideal of $\phi'|B$. For any maximal ideal $\Fq\not=\Fq_0$ of $B$ let $D_\Fq$ denote the commutant of $\End_{K^\sep}(\phi'|B)$ in $\End_{B_\Fq}(T_\Fq(\phi'|B))\cong\Mat_{r''\times r''}(B_\Fq)$, which is an order in a central simple algebra over $\Quot(B_\Fq)=E_\Fq$. Let $D_\Fq^1$ denote the multiplicative group of elements of $D_\Fq$ of reduced norm~$1$, which is a subgroup of $\SL_{r''}(B_\Fq)$. Choose an element $b_0\in B$ that generates a power of $\Fq_0$, view it as a scalar in $\prod_{\Fq\not=\Fq_0}\GL_{r''}(B_\Fq)$, and let $\overline{\langle b_0\rangle}$ denote the closure of the subgroup generated by it.
Let $K''\subset K^\sep$ be a finite extension of $K'$ over which all elements of $\End_{K^\sep}(\phi'|B)$ are defined.

\begin{Thm}\label{EndGalSpec3}
(Devic-Pink \cite[Thm.\;1.2]{DevicPink})
In the situation of Theorem \ref{EndGalSpec1}, the image of $\Gal(K^\sep/K'')$ in the adelic Galois representation associated to $\phi'|B$ is 
contained in $\prod_{\Fq\not=\Fq_0}D_\Fq^\times$ and commensurable with 
$$\overline{\langle b_0\rangle} \cdot \prod_{\Fq\not=\Fq_0}D_\Fq^1.$$
\end{Thm}

The images of Galois for $\phi'$ and $\phi$ up to commensurability can be determined from the image for $\phi'|B$ as explained in Devic-Pink \cite[\S6.2]{DevicPink}. Specifically, by \cite[Prop.\;6.7]{DevicPink} the characteristic ideal $\Fp_0'\subset A'$ of $\phi'$ is the unique maximal ideal of $A'$ above~$\Fq_0$. For each maximal ideal $\Fq\not=\Fq_0$ of $B$ there is a natural Galois equivariant isomorphism 
\UseTheoremCounterForNextEquation
\begin{equation}\label{TphiTphiB}
T_\Fq(\phi'|B) \ \cong\ \prod_{\Fp'|\Fq}T_{\Fp'}(\phi').
\end{equation}
% Since the center of $\End_{K^\sep}(\phi'|B)$ is~$B$, 
This induces a natural embedding 
\UseTheoremCounterForNextEquation
\begin{equation}\label{DFqTphi}
D_\Fq\ \into\ \prod_{\Fp'|\Fq}\End_{A'_{\Fp'}}(T_{\Fp'}(\phi')) 
\ \cong\ \prod_{\Fp'|\Fq}\Mat_{r'\times r'}(A'_{\Fp'}).
\end{equation}
Via Theorem \ref{EndGalSpec3} this determines the action of Galois on the Tate modules of~$\phi'$. A similar reduction process yields the action on the Tate modules of~$\phi$.

\medskip
We will use Theorem \ref{EndGalSpec2} to bound $E$ and $B$ from below, so the case (c') might cause us problems. But we can avoid these using the following additional result:

\begin{Prop}\label{Rank2}
In the situation of Theorem \ref{EndGalSpec1}, if $\rank(\phi')=2$, then $B=A'$.
\end{Prop}

\begin{Proof}
% To ease notation we replace $(A,K,\phi)$ by $(A',K',\phi')$. Then $\phi$ is non-isotrivial of special characteristic and of rank~$2$.
To ease notation we replace $K$ by~$K'$. Let $M$ denote the moduli scheme of Drinfeld $A'$-modules of rank~$2$, which is affine of relative dimension $1$ over~$\Spec(A')$. Since $\phi'$ is non-isotrivial, the associated $K$-valued point of $M$ lies over a generic point of the special fiber $M_{\Fp_0'}$ over $\Spec(A'/\Fp_0')$. 
On the one hand this shows that after replacing $K$ by a suitable subfield of $K^\sep$ we may assume that $K$ has transcendence degree~$1$ over~$\BF_p$. On the other hand, since $M_{\Fp_0'}$ is affine, 
% Since that fiber is an affine curve, 
there exists a place $v$ of~$K$ with local ring $\CO_v$ such that the $K$-valued point does not extend to a morphism $\Spec\CO_v\to M$. This means that $\phi'$ does not have potentially good reduction at~$v$. After replacing $K$ by a finite extension we may assume that $\phi'$ has semistable reduction at~$v$. Thus after conjugating $\phi'$ by an element of $K^\times$ we may assume that its coefficients are integral at $v$ and that its reduction has rank~$>0$.

Choose an extension of $v$ to $K^\sep$ and let $\hat\CO_v\subset K_v\subset K_v^\sep$ denote the corresponding completions of $\CO_v\subset K\subset K^\sep$. Let $I_v\subset D_v\subset\Gal(K^\sep/K)$ denote the respective inertia and decomposition groups.
Since $\phi'$ has rank~$2$, its Tate uniformization (see Drinfeld \cite[\S7]{Drinfeld1}) must consist of a Drinfeld $A'$-module $\psi_v$ of rank $1$ over $\Spec\hat\CO_v$ and an $A'$-lattice $\Lambda_v\subset K_v^\text{sep}$ of rank $1$ for the action of $A'$ on $K_v^\text{sep}$ via~$\psi_v$. Here by definition an $A'$-lattice is a finitely generated projective $A'$-submodule whose intersection with any ball of finite radius is finite. This implies that any non-zero element of $\Lambda_v$ has valuation $<0$. Also, since $\Lambda_v$ is finitely generated, after again replacing $K$ by a finite extension we may assume that $\Lambda_v\subset K_v$.

Take any maximal ideal $\Fp'\not=\Fp_0'$ of~$A'$. Then the Tate uniformization yields a natural $D_v$-equivariant short exact sequence
$$0 \longrightarrow T_{\Fp'}(\psi_v) \longrightarrow T_{\Fp'}(\phi') \longrightarrow \Lambda_v\otimes_{A'}A'_{\Fp'} \longrightarrow 0.$$
Here $I_v$ acts trivially on the outer terms; so in a suitable basis its action on $T_{\Fp'}(\phi')$ corresponds to a homomorphism
\UseTheoremCounterForNextEquation
\begin{equation}\label{inertunip}
I_v\ \longto\ U_{\Fp'}
\ :=\ \left(\!\begin{array}{cc}1&A'_{\Fp'}\\0&1\end{array}\!\!\right)
\ \subset\ \GL_2(A'_{\Fp'}).
\end{equation}
Let $\Delta$ denote the image of this homomorphism, viewed as a closed subgroup of the additive group of~$A'_{\Fp'}$.
%\begin{Lem}\label{Rank2Lem}
%\end{Lem}
%\begin{Proof}
%\end{Proof}
We claim that $\Delta$ is open in~$A'_{\Fp'}$.

To see this, we assume without loss of generality that the valuation $v$ is normalized on~$K_v$. Pick an element $\lambda\in\Lambda_v\setminus\{0\}$ and set $c:=-v(\lambda)\in\BZ^{\ge1}$. Recall that some power of $\Fp'$ is principal, say $\Fp^{\prime k}=(a')$ with $k>0$ and $a'\in A'$. The Tate uniformization thus yields a natural $D_v$-equivariant isomorphism 
$$\phi'[\Fp^{\prime k}]\ \cong\ 
\bigl\{ x \in K_v^\sep \bigm| 
        \psi_{v,a'}(x)\in\Lambda_v \bigr\} \bigm/ \Lambda_v.$$
Set $m:=\dim_{\BF_p}(A'/\Fp^{\prime k})$. Since $\psi_v$ is a Drinfeld $A'$-module of rank $1$ over $\hat\CO_v$, we have $\psi_{v,a'}=\sum_{i=0}^m u_i\tau^i$ with $u_i\in\hat\CO_v$ and $u_m\in\hat\CO_v^\times$. Therefore any solution $x\in K_v^\sep$ of the equation $\psi_{v,a'}(x)=\lambda$ satisfies $p^m\cdot v(x)=v(\lambda)=-c$. It follows that the field extension $K_v(x)/K_v$ has ramification degree at least $p^m/c$. The image of $I_v$ in the action on $\phi'[\Fp^{\prime k}]$ therefore also has order at least $p^m/c = |A'/\Fp^{\prime k}|/c$. 
But this image is naturally isomorphic to the image of $\Delta\subset A'_{\Fp'}$ in $A'/\Fp^{\prime k}$, which therefore has index at most~$c$.
Repeating the calculation with $\Fp^{\prime ki}$ in place of $\Fp^{\prime k}$ shows that for every integer $i>0$, the image of $\Delta\subset A'_{\Fp'}$ in $A'/\Fp^{\prime ki}$ has index at most~$c$. Passing to the inverse limit over $i$ we deduce that $\Delta\subset A'_{\Fp'}$ itself has index at most~$c$. It is therefore open, as claimed.

Now we can prove the proposition by contradiction. Suppose that $B\subsetneqq A'$, or equivalently $[F'/E]>1$. Then we can find a maximal ideal $\Fq\not=\Fq_0$ of $B$ and a maximal ideal $\Fp'$ of $A'$ above $\Fq$ such that $[F'_{\Fp'}/E_\Fq]>1$. We can also make $\Fq$ avoid the finitely many primes of~$E$ where the central simple $E$-algebra $\End^\circ_{K^\sep}(\phi'|B)$ is not split. Then its commutant $D_\Fq\otimes_{B_\Fq}E_\Fq$ is also split, i.e., isomorphic to the ring of $2\times2$-matrices over~$E_\Fq$. Theorem \ref{EndGalSpec3} with the embedding (\ref{DFqTphi}) thus implies that the image $\Gamma_\Fp$ of $\Gal(K^\sep/K)$ in the Galois representation on $T_{\Fp'}(\phi')$ is contained in a conjugate of $\GL_2(E_\Fq)$ in $\GL_2(F'_{\Fp'})$. But by the claim above $\Gamma_\Fp$ contains a conjugate of an open subgroup of $\binom{\kern2pt1\ \ A\rlap{$\scriptstyle{}'$}_{\Fp'}}{0\ \ \kern1pt 1\kern3pt}$. Together this is not possible with $[F'_{\Fp'}/E_\Fq]>1$, yielding the desired contradiction.
\end{Proof}

\medskip
{\bf Independence of Frobeniuses:} Next we will show that there exist Frobeniuses for $\phi'$ whose associated field extensions of $F'$ are maximally independent. This requires some group theoretical preparation. 

Consider a non\-archi\-me\-dean local field $L$ of equal characteristic $p$ with algebraic closure~$\oL$. Recall that an element of $\GL_{r'}(L)$ is called regular semisimple if it has $r'$ distinct eigenvalues in~$\oL$. Let us call an element totally split if its eigenvalues lie in~$L$, respectively totally inert if its eigenvalues generate an unramified field extension of degree $r'$ of~$L$. 

\begin{Lem}\label{GroupLem}
Every open subgroup of $\SL_{r'}(L)$ possesses an element $\gamma$ such that, for any $\delta\in \GL_{r'}(L)$ sufficiently close to~$\gamma$, every positive power of $\delta$ is regular semisimple and totally split. The same is true with totally inert in place of totally split.
\end{Lem}

\begin{Proof}
Let $\CO_L$ denote the valuation ring of $L$ and $(\pi)$ its maximal ideal. Choose $i\ge1$ such that the given subgroup contains all elements of $\SL_{r'}(\CO_L)$ which are congruent to the identity matrix modulo $(\pi^i)$. 
%If the characteristic of $L$ is zero but its residue characteristic is $p>0$, assume also that $\pi^i\in (p^2)$. Then any root of unity in $\oL$ which is congruent to $1$ modulo $(\pi^i)$ is trivial. 
Let $\gamma_0\in\GL_{r'}(\CO_L)$ be the diagonal matrix with diagonal entries $1+\pi^{i},1+\pi^{i+1},\ldots,1+\pi^{i+r'-1}$. Then $\gamma:=\gamma_0^{r'}\det(\gamma_0)^{-1}$ lies in the given subgroup of $\SL_{r'}(L)$. By construction $\gamma$ has $r'$ distinct eigenvalues in~$L$, which are all congruent to $1$ modulo $(\pi^i)$. For any $\delta\in \GL_{r'}(L)$ close to~$\gamma$, the characteristic polynomial of $\delta$ is close to that of~$\gamma$. But by Hensel's lemma split separable polynomials remain split separable under small deformations. Thus any $\delta\in \GL_{r'}(L)$ sufficiently close to $\gamma$ has $r'$ distinct eigenvalues in~$L$, which are all congruent to $1$ modulo $(\pi^i)$. Moreover, if some positive power $\delta^n$ had two equal eigenvalues, two eigenvalues of $\delta$ would differ by a nontrivial root of unity congruent to $1$ modulo $(\pi)$, which does not exist. Thus $\delta^n$ is regular semisimple and totally split, as desired.

To prove the same assertion with totally inert in place of totally split, let $L'$ be an unramified extension of degree $r'$ of~$L$ with valuation ring $\CO_{L'}$. Choose an element $\alpha\in\CO_{L'}$ whose residue class generates the residue field extension $k'/k$ and has trace $\trace_{k'/k}(\alpha)=0$.
% with $\CO_{L'}=\CO[\alpha]$ and $\trace_{L'/L}(\alpha)=0$. 
Identify $\CO_{L'}$ with a subring of the matrix ring $\Mat_{r'\times r'}(\CO_L)$, and set $\gamma_0:={1+\pi^i\alpha}\in\GL_{r'}(\CO_L)$. Then $\det(\gamma_0)\equiv1+\pi^i\trace_{k'/k}(\alpha)\equiv1$ modulo $(\pi^{i+1})$. Dividing one matrix column of $\gamma_0$ by this determinant yields an element $\gamma\in\SL_{r'}(\CO_L)$ which is congruent to $1+\pi^i\alpha$ modulo $(\pi^{i+1})$. Thus $\gamma$ lies in the given subgroup. Consider any $\delta\in\GL_{r'}(\CO_L)$ congruent to $\gamma$ modulo $(\pi^{i+1})$. Then $(\delta-1)/\pi^i$ has coefficients in $\CO_L$ and is congruent to $\alpha$ modulo $(\pi)$; hence its residue class generates $k'$ over~$k$. Thus the $\CO_L$-subalgebra of $\Mat_{r'\times r'}(\CO_L)$ generated by it is isomorphic to $\CO_{L'}$. It follows that the $L$-subalgebra of $\Mat_{r'\times r'}(L)$ generated by $\delta$ is isomorphic to~$L'$; hence $\delta$ is regular semisimple and totally inert. Moreover, the ratio of any two distinct eigenvalues of $\delta$ is congruent to $1$ modulo $(\pi)$ and therefore not a root of unity. Thus any positive power $\delta^n$ is again regular semisimple and generates the same $L$-subalgebra, hence is again totally inert, as desired.
\end{Proof}

\medskip
Now we return to Drinfeld modules of arbitrary characteristic, keeping the notation from before. For any maximal ideal $\Fm'$ of $R'$ we abbreviate $F'_{\Fm'}:=\End_{k_{\Fm'}^\sep}^\circ(\phi'_{\Fm'})$.

\begin{Prop}\label{X1}
There exist maximal ideals $\Fm',\Fn'\subset R'$, such that $F'_{\Fm'}$ and $F'_{\Fn'}$ are commutative and linearly disjoint over~$F'$, that is, their tensor product over $F'$ is a field.
\end{Prop}

\begin{Proof}
% As before set $r':=\rank(\phi')$. 
If $r'=1$, then $F'_{\Fm'}=F'$ for all~$\Fm'$ and the assertion is trivial. So assume that $r'>1$. Then $\phi'$ is not isotrivial. If $\phi'$ has special characteristic, let $B$, $E$, and $D_\Fq$ be as above. Otherwise, set $B:=A'$ and $E:= F'$ and $D_\Fq :=\Mat_{r'\times r'}(B_\Fq)$. 

Let $F''\subset F'$ be the maximal subfield which is separable over~$E$. Then there exist infinitely many maximal ideals $\Fq$ of $B$ which are totally split in~$F''$. For almost all of these we also have $\Fq\not=\Fq_0$ and $D_\Fq \cong\Mat_{r'\times r'}(B_\Fq)$. We select two distinct maximal ideals $\Fq$ and $\Fq'$ of $B$ with these properties. Let 
$$\tilde\Gamma\ \subset\ D_\Fq^\times\times D_{\Fq'}^\times \ \cong\ 
\GL_{r'}(B_\Fq)\times\GL_{r'}(B_{\Fq'})$$
denote the image of $\Gal(K^\sep/K'')$ in the Galois representation on $T_\Fq(\phi'|B)\times T_{\Fq'}(\phi'|B)$. Then Theorem \ref{EndGalGen}, respectively \ref{EndGalSpec3}, implies that $\tilde\Gamma$ contains an open subgroup of $\SL_{r'}(B_\Fq)\times\SL_{r'}(B_{\Fq'})$.

\begin{Lem}\label{X1Lem}
There exists a maximal ideal $\Fm'$ of $R'$ such that $F'_{\Fm'}$ is commutative and any maximal ideal of $A'$ above $\Fq$ is totally split in $F'_{\Fm'}$, while any maximal ideal of $A'$ above $\Fq'$ is totally inert in~$F'_{\Fm'}$.
\end{Lem}

\begin{Proof}
Using Lemma \ref{GroupLem} choose an element $\gamma\in\SL_{r'}(B_\Fq)$ close to the identity element, such that for any $\delta\in \GL_{r'}(B_\Fq)$ sufficiently close to~$\gamma$, every positive power of~$\delta$ is regular semisimple and totally split. Likewise choose an element $\gamma'\in\SL_{r'}(B_{\Fq'})$ close to the identity element, such that for any $\delta'\in \GL_{r'}(B_{\Fq'})$ sufficiently close to~$\gamma'$, every positive power of~$\delta'$ is regular semisimple and totally inert. As these elements can be chosen arbitrarily close to the identity element, we can require that $\tilde\gamma := (\gamma,\gamma')$ is an element of~$\tilde\Gamma$. 
Since the images of Frobenius elements in $\Gal(K^\sep/K'')$ form a dense subset of~$\tilde\Gamma$, there then exists a maximal ideal $\Fm'$ of $R'$ such that the image $\tilde\delta=(\delta,\delta')$ of $\Frob_{\Fm'}$ satisfies the stated conditions, i.e., any positive power of $\delta$ is regular semisimple and totally spilt and any positive power of $\delta'$ is regular semisimple and totally inert.
We claim that $\Fm'$ has the desired properties.

To see this recall that $\End_{k_{\Fm'}^\sep}^\circ(\phi'_{\Fm'}|B)=\End_{\ell_{\Fm'}}^\circ(\phi'_{\Fm'}|B)$ for some finite field extension $\ell_{\Fm'}\subset k_{\Fm'}^\sep$ of $k_{\Fm'}$, say of degree $n\ge1$. 
Its center $E_{\Fm'}$ is thus the field extension of $E$ which is generated by $\Frob_{\Fm'}^n$. 
Moreover, the minimal polynomial of $\Frob_{\Fm'}^n$ over $E$ is equal to that of $\delta^n$ and of~$\delta^{\prime n}$. As these elements are regular semisimple, it follows that $\Frob_{\Fm'}^n$ is separable of degree $r'$ over~$E$. Thus $E_{\Fm'}$ is separable of degree $r'$ over~$E$. 

Next, the reduction of endomorphisms (\ref{EndRed}) induces a natural homomorphism of $E_{\Fm'}$-algebras
\UseTheoremCounterForNextEquation
\begin{equation}\label{X1Red}
E_{\Fm'}\otimes_E\End^\circ_{K^\sep}(\phi'|B) \longto \End_{k_{\Fm'}^\sep}^\circ(\phi'_{\Fm'}|B).
\end{equation}
Recall that $r'=\rank(\phi')$, so that $r''=\rank(\phi'|B) = r'd$ with $d:=[F'/E]$. Then $\End^\circ_{K^\sep}(\phi'|B)$ is a central simple $E$-algebra of dimension~$d^2$. 
On the other hand, since $E_{\Fm'}$ is the center of $\End_{k_{\Fm'}^\sep}^\circ(\phi'_{\Fm'}|B)$ and of degree $r'$ over~$E$, the equation $\rank(\phi'|B) = r''= r'd$ implies that $\dim_{E_{\Fm'}}(\End_{k_{\Fm'}^\sep}^\circ(\phi'_{\Fm'}|B))\le d^2$. Thus the source and target in (\ref{X1Red}) are central simple $E_{\Fm'}$-algebras of dimension~$d^2$, respectively $\le d^2$; hence the homomorphism is an isomorphism.

Now observe that by the definition of endomorphisms $\End^\circ_{K^\sep}(\phi')$ is simply the commutant of $F'$ within $\End^\circ_{K^\sep}(\phi'|B)$. The fact that $\End^\circ_{K^\sep}(\phi')=F'$ thus means that $F'$ is a maximal commutative subalgebra of~$\End^\circ_{K^\sep}(\phi'|B)$. Therefore the isomorphism (\ref{X1Red}) maps $E_{\Fm'}\otimes_EF'$ isomorphically to a maximal commutative subalgebra of $\End_{k_{\Fm'}^\sep}^\circ(\phi'_{\Fm'}|B)$. But again by the definition of endomorphisms $\End_{k_{\Fm'}^\sep}^\circ(\phi'_{\Fm'})$ is simply the commutant of $F'$ within $\End_{k_{\Fm'}^\sep}^\circ(\phi'_{\Fm'}|B)$. 
As the center of $\End_{k_{\Fm'}^\sep}^\circ(\phi'_{\Fm'}|B)$ is~$E_{\Fm'}$, this commutant is equal to the commutant of the image of $E_{\Fm'}\otimes_EF'$, and hence equal to the image of $E_{\Fm'}\otimes_EF'$. This shows that $F'_{\Fm'}:=\End_{k_{\Fm'}^\sep}^\circ(\phi'_{\Fm'})$ is isomorphic to $E_{\Fm'}\otimes_EF'$ over~$F'$. In particular $F'_{\Fm'}$ is commutative.

Finally, the fact that $\delta^n$ is totally split implies that $\Fq$ is totally split in the field extension $E_{\Fm'} = E(\Frob_{\Fm'}^n)$. It follows that any maximal ideal of $A'$ above $\Fq$ is totally split in ${E_{\Fm'}\otimes_EF'}\allowbreak\cong F'_{\Fm'}$. 
Likewise, the fact that $\delta^{\prime n}$ is totally inert implies that $\Fq'$ is totally inert in~$E_{\Fm'}$. 
Since by assumption $\Fq'$ is totally split in the maximal separable subextension of $F'/E$, every maximal ideal of $A'$ above $\Fq'$ has the same residue field as~$\Fq$. Thus every maximal ideal of $A'$ above $\Fq'$ is totally inert in $E_{\Fm'}\otimes_EF'\cong F'_{\Fm'}$. Therefore ${\Fm'}$ has all the desired properties.
\end{Proof}

\medskip
To finish the proof of Proposition \ref{X1}, choose any $\Fm'$ as in Lemma \ref{X1Lem}. Applying Lemma \ref{X1Lem} with the roles of $\Fq$ and $\Fq'$ reversed, we also choose a maximal ideal $\Fn'$ of $R$ such that $F'_{\Fn'}$ is commutative and that any maximal ideal of $A'$ above $\Fq'$ is totally split in~$F'_{\Fn'}$, while any maximal ideal of $A'$ above $\Fq$ is totally inert in~$F'_{\Fn'}$. Together these properties imply that $F'_{\Fm'}$ and $F'_{\Fn'}$ are linearly disjoint over~$F'$, and we are done.
\end{Proof}

%%%%%%%%%%%%%%%%%%%%%%%%%%%%%%%%%%%%%%%%%%%%%%%%%%%%%%%%%%%%%%%%%%%%%%%%%%%
%\newpage
\section{Computer algebra prerequisites}
\label{Comp}

In this section we briefly recall the methods from computer algebra which are used in the rest of the article. As a general reference, one can consult for example the book ``Computational Commutative Algebra 1'' by Kreuzer and Robbiano \cite{KreuzerRobbiano}. Many of the operations mentioned here are implemented in common computer algebra systems.
%(for instance magma or macaulay2).  

\medskip{\bf Representation of algebras and fields:}
Any finitely generated $\BF_p$-algebra $R$ can be represented as the quotient of a polynomial ring $\BF_p[\UX] := \BF_p[X_1,\ldots,X_r]$ by a finitely generated ideal~$J$. 
Using Gr\"obner bases one can effectively decide whether $J$ is prime, or equivalently whether $R$ is integral. 
The localization of $R$ with respect to finitely many elements $x_1,\ldots,x_s$ can be represented on the same footing as $R':=R[1/x_1\cdots x_s]=R[Y]/(x_1\cdots x_sY-1)$.

Any finitely generated field $K$ over $\BF_p$ can be represented as the field of fractions of $R:=\BF_p[\UX]/J$ for a prime ideal~$J$. Any calculation with ideals in $K[Y_1,\ldots,Y_s]$ reduces to one in $R'[Y_1,\ldots,Y_s]$ for a suitable localization $R'$ of~$R$.

%Using Gr\"obner bases one can effectively compute the radical and the primary decomposition of~$J$ and decide whether $J$ is prime. 

\medskip{\bf Basic operations on elements:} 
Let $R=\BF_p[\UX]/J$ be a finitely generated $\BF_p$-algebra. Using a Gr\"obner basis of $J$, for every element of $R$ one can compute its unique reduced representative with respect to this basis. Thus one can effectively decide whether two given elements of $R$ are equal. If $R$ is integral, one can therefore also decide whether two elements of its field of fractions are equal. 

Using Gr\"obner bases one can also test whether a given element is contained in a given ideal of~$R$. In particular, if $R$ is integral, one can test whether one element divides another in $R$ and, if so, determine the quotient. Thus one can decide whether an element of the field of fractions already lies in~$R$.

\medskip{\bf Ideals and subrings:}
For any homomorphism $f\colon S \to R$ of finitely generated $\BF_p$-algebras and any ideal $J$ of $R$ one can effectively determine the ideal $f^{-1}(J)$ of~$S$. In particular, one can determine $\Ker(f)$ and hence obtain an explicit representation of $\Image(f)\cong S/\Ker(f)$. Applying this when $S$ is a polynomial ring over~$\BF_p$, one can thus explicitly describe the subalgebra generated by finitely many given elements of~$R$.

Furthermore, one can construct a sequence of all maximal ideals of $R$. 

\medskip{\bf Normalization:}
If $R$ is integral with field of fractions~$K$, one can effectively describe the normalization $R'$ of $R$ together with the inclusions $R\into R'\into K$, see for example Singh-Swanson  \cite{SinghSwanson}. 

\medskip{\bf Field extensions:}
Let $K$ be a field which is finitely generated over~$\BF_p$. Then for any irreducible polynomial $P\in K[X]$ one can write down a field extension of $K$ generated by a root of~$P$, namely as $K[X]/(P)$. 
Given an arbitrary polynomial $P\in K[X]$, one can effectively find its irreducible factors with multiplicities by Steel \cite{Steel}. By iteration one can therefore effectively describe a splitting field of $P$ over~$K$.
% By looking at the irreducible factors of $P$ one can decide whether that extension is separable.

For any field extension $K\subset L$ and any element $x\in L$, one can effectively decide whether $x$ is algebraic over $K$ and, if so, determine its minimal polynomial over~$K$. By factoring a polynomial over $L$ one can determine all its roots in~$L$. In particular, one can therefore determine all conjugates of $x$ over~$K$ in~$L$.

Also, for any simple finite extension $K\subset K'$, one can effectively describe all homomorphisms $K'\to L$ over~$K$, by mapping the generator of $K'$ to roots of its minimal polynomial. By iteration over simple extensions, one can effectively describe all homomorphisms $K'\to L$ over~$K$ for any finite extension $K\subset K'$.

Moreover, one can effectively decide whether two field extensions $K'/K$ and $K''/K$ are linearly disjoint in that their tensor product $K'\otimes_KK''$ is a field. 
Indeed, if $K'/K$ is simple, this is equivalent to the minimal polynomial of the generator over $K$ remaining irreducible over~$K''$. The general case follows by iteration over simple extensions.
%Indeed, if each $K_i$ is given as the field of fractions of a finitely generated integral $\BF_p$-algebra $R_i$, one can effectively compute the kernel of the natural homomorphism $R_1\otimes_{R_0}R_2\to K_1\otimes_{K_0}K_2$ and decide whether it is a prime ideal.

\medskip
% One can determine the field of constants of~$K$.

For any field $K$ we let $K^\sep\subset\overline{K}$ denote a separable, respectively an algebraic closure of~$K$. Though one cannot effectively construct these and compute in them, one can calculate in any finite extension and enlarge it whenever necessary. 
Throughout, all finite separable extensions of $K$ are tacitly assumed to be contained in~$K^\sep$.

\medskip{\bf Solving polynomial equations:}
Let $\CS$ be a system of finitely many polynomial equations in several variables over $K$ which is known to have only finitely many solutions in~$\oK$. Then one can determine a finite extension $K'$ of $K$ such that all solutions of $\CS$ lie in $K'$ and one can find those solutions; Lazard \cite{Lazard} gives a possible way of doing this. 

\medskip{\bf Intermediate fields:}
For any finite separable field extension $K\subset L$ one can effectively find a Galois closure~$\tilde L$ and determine the Galois group of $\tilde L/K$. For every subgroup of this Galois group one can effectively determine generators of the associated intermediate field. In this way one can make a finite list of all intermediate fields of $L/K$. 

More generally let $L/K$ be an arbitrary finite field extension with maximal separable subextension $L'/K$. Then any intermediate field of $L/K$ is a purely inseparable field extension of an intermediate field of $L'/K$. 
If $K$ has transcendence degree $1$ over~$\BF_p$, any purely inseparable extension is determined by its degree and generated by $p$-power roots; hence one can also make a finite list of all intermediate fields of $L/K$ in this case.

\medskip{\bf Transcendence degree 1:} 
We will often deal with finitely generated integral domains over $\BF_p$ of transcendence degree~$1$. Any such ring $B$ possesses a transcendent element $t$ such that $B$ is a finitely generated $\BF_p[t]$-module. One can thus present $B$ efficiently via a basis as $\BF_p[t]$-module and a multiplication table with entries in $\BF_p[t]$.
%For any other transcendent element $t'\in B$ such that $B$ is a finitely generated $\BF_p[t']$-module, 
For any other such element $t'\in B$ one can translate this presentation over $\BF_p[t]$ into one over $\BF_p[t']$, using commutative algebra over~$\BF_p[t,t']$. 

In the same way one can describe any torsion free commutative or non-commutative $B$-algebra which is finitely generated as a $B$-module. This reduces many computations with modules and ideals to linear algebra over the principal ideal domain $\BF_p[t]$. 

\medskip{\bf Modules:} 
For any finitely generated $B$-module $M$ one can find its rank and the elementary divisors as an $\BF_p[t]$-module. In particular, one can decide whether $M$ is finite and, if so, make a list of its elements. If $M$ is torsion free, for any submodule $N$ of $M$ one can effectively compute the saturation
$\{m\in M\mid \exists\,b\in B\setminus\{0\}\colon bm\in N\}$. 
In particular, for any $\Quot(B)$-subspace $V$ of $M\otimes_B\Quot(B)$ one can determine $V\cap M$.

\medskip{\bf Admissible coefficient rings:}
Let $F$ be the function field of an irreducible smooth projective curve $C$ over~$\BF_p$. Let $A$ be the subring of elements of $F$ which are regular outside a fixed point $\infty$ of~$C$. We call such $A$ an admissible coefficient ring.
For any ideal $\Fa\subset A$, one can compute its prime factorization and the number $\dim_{\BF_p}(A/\Fa)$. One can also compute the degree  $\deg(\infty)$ of the residue field at $\infty$ over~$\BF_p$. 
For any integer $n\geq 0$, the finite set of elements $a$ of $A$ with $a=0$ or $\dim_{\BF_p}(A/Aa)\leq d$ is just the Riemann-Roch space $\Gamma\bigl(C,\CO_C(\lfloor\frac{n}{\deg(\infty)}\rfloor\infty)\bigr)$, which can be effectively determined by Hess \cite{Hess}.

%\begin{Lem}\label{RootofUnityLem}
%For any monic polynomial $P$ over~$F$, one can effectively decide whether, writing $P(X)=\prod_{i=1}^r(X-\alpha_i)$ over an algebraic closure $\oF$ of~$F$, for all $n\ge1$ the elements $\alpha_1^n,\ldots,\alpha_r^n$ are pairwise distinct.
%\end{Lem}
%
%\begin{Proof}
%If $X^2$ divides~$P(X)$, the answer is no. Otherwise, if $X$ divides~$P(X)$, we can replace $P(X)$ by $P(X)/X$ without changing the answer. Thus we may assume that $P(0)\not=0$ and hence all $\alpha_i\not=0$. The condition is then equivalent to saying that for all $i\not=j$ the ratio $\alpha_i/\alpha_j$ is not a root of unity. So consider the polynomial
%$$Q(X)\ :=\prod_{i\not=j}\Bigl(X-\frac{\alpha_i}{\alpha_j}\Bigr).$$
%Being symmetric functions in $\alpha_1,\ldots,\alpha_r$, the coefficients of $Q$ lie in $F$ and can be effectively computed from those of~$P$. The condition is equivalent to saying that no root of $Q$ is algebraic over~$\BF_p$. So if $\BF$ denotes the constant field of~$F$, the condition is equivalent to saying that the factorization of $Q$ into monic irreducible polynomials over $F$ possesses no factor defined over~$\BF$. But by ... this factorization can be effectively computed, and by ... it can be effectively decided whether one factor is defined over~$\BF$.
%\end{Proof}

\medskip
We finish this section with two more specialized facts.

\begin{Lem}\label{FieldExtLem}
For any simple finite field extension $F(x)$ of~$F$, there exists an integer $m\ge1$ such that $\bigcap_{n\ge1}F(x^n)=F(x^m)$. Moreover, knowing the minimal polynomial of $x$ over~$F$ one can effectively find such $m$ as well as the minimal polynomial and the degree of $x^m$ over~$F$.
In particular, one can effectively decide whether $F(x^n)=F(x)$ for all $n\ge1$.
\end{Lem}

\begin{Proof}
Choose $m$ such that $[F(x^m)/F]$ is minimal. Then for any $n\ge1$, we have 
%$F(x^{nm})\subset F(x^m)$ and $[F(x^{nm})/F]=[F(x^m)/F]$ and hence $F(x^{nm})=F(x^m)$. Therefore 
$F(x^m)=F(x^{nm})\subset F(x^n)$. This proves the first statement of the lemma.

To find $m$ effectively, note first that if $x=0$, then $m=1$ does the job. Otherwise, let $P(X)$ be the minimal polynomial of $x$ over~$F$. By looking at the coefficients of $P(X)$ one can find the largest power $p^i$ such that $P(X)=Q(X^{p^i})$ for some polynomial $Q(X)$. Then $Q(X)$ is the minimal polynomial of $x^{p^i}$ over~$F$. After replacing $x$ by $x^{p^i}$ we can thus assume that $P(X)\not\in F[X^p]$, in other words that $P(X)$ is separable. 
Let $\sigma_1,\ldots,\sigma_r$ be the pairwise distinct homomorphisms  over $F$ from $F(x)$ into a separable closure $F^\sep$ of~$F$. Then $P(X)=\prod_{i=1}^r(X-\sigma_i(x))$. Consider the polynomial
$$Q(X)\ :=\prod_{i\not=j}\Bigl(X-\frac{\sigma_i(x)}{\sigma_j(x)}\Bigr).$$
Being symmetric functions in $\sigma_1(x),\ldots,\sigma_r(x)$, the coefficients of $Q$ lie in $F$ and can be effectively computed from those of~$P$. We can then effectively compute the factorization of $Q(X)$ into monic irreducible polynomials over~$F$. By determining which of their coefficients are algebraic over $\BF_p$ we can effectively decide which of these factors are defined over the constant field of~$F$. From these we can effectively find a positive integer $m$ such that all their roots are $m$-th roots of unity. This integer has the property that for all $n\ge1$ and $i\not=j$, if $\sigma_i(x^n)=\sigma_j(x^n)$, then $\sigma_i(x^m)=\sigma_j(x^m)$. By Galois theory this implies that $F(x^m)\subset F(x^n)$, so $m$ has the desired property.

Using symmetric functions again one can now effectively compute the polynomial $\prod_{i=1}^r(X-\sigma_i(x^m))$ in $F[X]$. As this is a power of the minimal polynomial of $x^m$ over~$F$, by factorization one can effectively determine this minimal polynomial, and hence also its degree. This proves the second statement of the lemma. 

In particular, knowing the degrees of of the minimal polynomials of $x$ and $x^m$ over $F$ one can effectively decide whether $F(x^m)=F(x)$. This implies the last statement.
\end{Proof}

\medskip
{\bf Computation in $K[\tau]$:}
As before let $K$ be a finitely generated field over~$\BF_p$. By definition an element $u=\sum_{i=0}^nu_i\tau^i\in K[\tau]$ with $u_n\not=0$ has degree $\deg_\tau(u):= n$, and the zero element has degree $-\infty$. This degree is additive in products, that is, for any $u,v\in K[\tau]$ we have $\deg_\tau(uv)=\deg_\tau(u)+\deg_\tau(v)$. Also, any left ideal of the ring $K[\tau]$ is principal, and $K[\tau]^\times=K^\times$.
% Goss \cite[\S1.6]{GossBS}

\begin{Prop}\label{gcrd&lclm}
For any elements $u,v\in K[\tau]$ with $v\not=0$ there exist unique $q,r\in K[\tau]$ with $u=qv+r$ and $\deg_\tau(r)<\deg_\tau(v)$.
Any finite subset of $K[\tau]$ possesses a greatest common right divisor and a least common left multiple, which are unique up to left multiplication by an element of~$K^\times$.
All of these can be computed effectively.
\end{Prop}

\begin{Proof}
For the first statement, that $K[\tau]$ is euclidean with respect to right division, see Goss \cite[\S1.6]{GossBS}. That $q$ and $r$ can be computed effectively is shown as in a commutative polynomial ring, for instance by induction on $\deg_\tau(u)$ and comparison of the highest coefficients.
% because $\deg_\tau(f)\le\deg_\tau(u)-\deg_\tau(v)$ and 
The corresponding euclidean algorithm yields the greatest common right divisor of any two, and consequently of any finite number of, elements of $K[\tau]$. 

For the least common left multiple of a finite subset $\CS\subset K[\tau]$ consider the left $K[\tau]$-module $M := \bigoplus_{u\in\CS}K[\tau]/K[\tau]u$.
The common left multiples of $\CS$ are precisely those elements of $K[\tau]$ which annihilate the element $(1+K[\tau]u)_{u\in\CS}$ of~$M$. Thus they form a left ideal of $K[\tau]$, which is therefore principal. Any generator of this ideal is a least common left multiple of~$\CS$. If $0\in \CS$, this least common left multiple is $0$; otherwise its degree in $\tau$ is at most $\dim_K(M)=\sum_{u\in \CS}\deg_\tau(u)<\infty$. In this case it can be determined by a bounded number of polynomial calculations over~$K$.
\end{Proof}

%%%%%%%%%%%%%%%%%%%%%%%%%%%%%%%%%%%%%%%%%%%%%%%%%%%%%%%%%%%%%%%%%%%%%%%%%%%
%\newpage
\section{Bits of algorithms}
\label{AlgBits}

Throughout this section we fix a Drinfeld $A$-module $\phi\colon A\to K[\tau]$ over a finitely generated field~$K$.

\begin{Prop}\label{CharIsoEff}
One can effectively determine the rank and height and characteristic ideal $\Fp_0$ of~$\phi$, hence whether $\phi$ has generic or special characteristic, and whether $\phi$ is ordinary resp.\ isotrivial.
\end{Prop}

\begin{Proof}
Choose a non-constant $t\in A$ and write $\phi_t=\sum_{i=0}^{n}x_i\tau^i$ with $x_n\not=0$. By definition the rank of $\phi$ is the quotient $n/\dim_{\BF_p}(A/At)$. By our computer algebra prerequisites, it can therefore be determined effectively. 
Also $\phi$ has generic characteristic if and only if $x_0$ is transcendental over~$\BF_p$. Specifically $\Fp_0$ is defined as the kernel of the homomorphism $A\to K$, $a\mapsto d\phi_a$. By our computer algebra prerequisites, it can be determined effectively. 

If $\Fp_0\not=0$, one can choose a new element $t\in\Fp_0\setminus\{0\}$. Write $\phi_t=\sum_{i=m}^{n}x_i\tau^i$ with $x_m\not=0$, and write $At=\Fp_0^k\Fa$ for an ideal $\Fa$ prime to $\Fp_0$. Then by definition the height of $\phi$ is the quotient $m/\dim_{\BF_p}(A/\Fp_0^k)$, and $\phi$ is ordinary if and only if its height is~$1$. This can therefore also be determined effectively. 

%Choose a non-constant $b\in A$ and write $\phi_b=\sum_{i=m}^{n}x_i\tau^i$ with $x_m\not=0$ and $x_n\not=0$. By our computer algebra prerequisites, the number $d:=\dim_{\BF_p}(A/Ab)$ can be determined effectively. By definition $\phi$ has rank $n/d$ and height $m/d$, and $\phi$ is ordinary if and only if the height is~$1$; hence all of these can be determined effectively. Also $\phi$ has generic characteristic if and only if $x_0$ is transcendental over~$\BF_p$. More specifically, by definition $\Fp_0$ is the kernel of the homomorphism $A\to K$, $a\mapsto d\phi_a$. By our computer algebra prerequisites, this can be determined effectively.

Next $\phi$ is isotrivial if and only if there exists $y\in\oK^\times$ such that $y^{-1}\phi_ay$ has coefficients in a finite field for every $a\in A$. We claim that it is enough to check this condition for the chosen element $t\in A$. Indeed, if it holds for~$t$, after replacing $\phi$ by $y^{-1}\phi y$ we may assume that $\phi_t$ has coefficients in a finite field $k\subset\oK$. Setting $B':=\BF_p[t]$, the restriction $\phi|B'$ is then a Drinfeld $B'$-module defined over~$k$. Thus there exists a finite extension $k'\subset\oK$ of $k$ with $\End_\oK(\phi|B')=\End_{k'}(\phi|B')\subset k'[\tau]$. Since $\phi$ induces an embedding $A\into\End_\oK(\phi|B')$, it follows that $\phi$ itself is defined over~$k'$ and is therefore isotrivial, as claimed.

To test the condition for~$t$, observe that $y^{-1}\phi_ty=\sum_{i=0}^{n}x_iy^{q^i-1}\tau^i$. A direct calculation shows that there exists $y\in\oK{}^\times$ such that all $x_iy^{q^i-1}$ are algebraic over~$\BF_p$ if and only if the ratios $x_i^{q^n-1}/x_n^{q^i-1}$ are algebraic over $\BF_p$ for all $0\le i<n$. By our computer algebra prerequisites, this condition can be tested effectively.
\end{Proof}

\begin{Prop}\label{IsogEff}
One can effectively find $A'$ and $\phi'$ and $h$ in Proposition \ref{Isog}, as well as an isogeny $h'\colon \phi'|A\to\phi$ over $K$ and an element $a\in A\setminus\{0\}$ with $h'h=\phi_a$.
\end{Prop}

\begin{Proof}
We follow the construction in Devic-Pink \cite[Prop.\;4.3]{DevicPink}.
By our computer algebra prerequisites, one can effectively describe the normalization $A'$ and, since $A'/A^\smalldash$ is finite, find an element $a\in A\setminus\{0\}$ satisfying $A'a\subset A^\smalldash$. Choose a system $\CR$ of non-zero representatives of cosets of $A'$ modulo~$A^\smalldash$. 

For simplicity we denote the given embedding $A^\smalldash\into\End_K(\phi)$ by $a^\smalldash\mapsto\phi_{a^\smalldash}$. 
Using Proposition \ref{gcrd&lclm}, for any $a'\in\CR$ one can effectively find the least common left multiple of $\phi_a$ and $\phi_{a'a}$ in $K[\tau]$ and write it in the form $h_{a'}\phi_{a'a}$ for some $h_{a'}\in K[\tau]$. 
Since $A'$ is commutative, we have $\phi_a\phi_{a'a}=\phi_{a'a}\phi_a$. Thus this element is a common left multiple of $\phi_{a'a}$ and $\phi_a$, hence a left multiple of $h_{a'}\phi_{a'a}$; hence $\phi_a$ is a left multiple of~$h_{a'}$. Next, using Proposition \ref{gcrd&lclm} again one can effectively find the least common left multiple $h\in K[\tau]$ of the elements $h_{a'}$ for all $a'\in\CR$. 
Since $\phi_a$ is already a common left multiple of these elements, it is a left multiple of~$h$; in other words we have $\phi_a=h'h$ for some $h'\in K[\tau]$. 

This construction has the following effect. For any non-zero element $u\in K[\tau]$ consider the finite subgroup scheme $\ker(u)$ of $\BG_{a,K}$. For any two non-zero elements $u,v\in K[\tau]$ with least common left multiple $wv$ we then have $\ker(u)+\ker(v)=\ker(wv)$ and hence 
$$v(\ker(u))\ =\ v(\ker(u)+\ker(v))\ =\ v(\ker(wv))\ =\ \ker(w).$$ 
In the above construction we therefore have $\phi_{a'a}(\ker(\phi_a))=\ker(h_{a'})$. Summing over all $a'\in\CR$ implies that 
$$\sum_{a'\in\CR}\phi_{a'a}(\ker(\phi_a)) \ =\ \ker(h).$$
Thus $h$ is the element $f$ from the proof of \cite[Prop.\;4.3]{DevicPink}. There we constructed a unique Drinfeld $A'$-module $\phi'\colon A'\to K[\tau]$ such that $h$ is an isogeny $\phi\to\phi'|A$. Since $h'h=\phi_a$ with $a\not=0$, it follows that $h'$ is an isogeny $\phi'|A\to\phi$. 
To find $\phi'$ explicitly, for any element $a'\in A'$ we have $a'a\in A^\smalldash$ and can therefore calculate 
$$\phi'_{a'}h\phi_a = \phi'_{a'}\phi'_ah  = \phi'_{a'a}h = h\phi_{a'a}.$$
This means that $\phi'_{a'}$ is the quotient of $h\phi_{a'a}$ upon right division by $h\phi_a$ in $K[\tau]$, which can again be computed explicitly by Proposition \ref{gcrd&lclm}.
\end{Proof}

\begin{Prop}\label{MinCoeff}
Let $f\in\End_K(\phi)$. As an element of the commutative subfield $F(f)$ of $\End^\circ_K(\phi)$ of finite degree over~$F$, it possesses a unique minimal polynomial $\min_f\in F[X]$. This polynomial has coefficients in~$A$, and writing $\min_f(X) = X^k+a_1X^{k-1}+\ldots+a_k$, for each $1\le i\le k$ we have 
$$\deg_\tau(\phi_{a_i})\le i\cdot\deg_\tau(f),$$
with equality for $i=k$.
\end{Prop}

\begin{Proof}
Since $\End_K(\varphi)$ is finitely generated as an $A$-module and $A$ is integrally closed, the minimal polynomial has coefficients in $A$.

Let $A'$ be the integral closure of $A[f]$ within $F':=F[f]$. By Proposition \ref{Isog}, this is an admissible coefficient ring and there exist a Drinfeld $A'$-module $\phi'$ and an isogeny $h\colon\phi\to\phi'|A$ over~$K$. This induces an isomorphism of $F$-algebras $\epsilon\colon\End_K^\circ(\phi)\stackrel{\sim}{\to}\End_K^\circ(\phi'|A)$ satisfying $\epsilon(g)h=hg$ for all $g\in\End_K^\circ(\phi)$; in particular $\epsilon(f)h=hf$. Thus $f$ and $\epsilon(f)=\phi'_f$ have the same degree in $\tau$ and the same minimal polynomial over~$F$. After replacing $\phi$ by $\phi'|A $, we may therefore assume that $A'\subset \End_K(\varphi)$.

Since $A\subset A'$ is an inclusion of admissible coefficient rings, there is a unique place $\infty'$ of $F'$ lying over $\infty$. For $f\in F'$ this implies that $\min_f$ has a unique slope at~$\infty$. More precisely, let $\ord_{\infty'}$ be the normalized valuation on $F'$ associated to~$\infty'$ and $\deg(\infty')$ the degree of its residue field over~$\BF_p$. Then by the definition of $\rank(\phi')$ we have 
$$\deg_\tau(\phi'_{a'}) 
\ =\ \rank(\phi')\cdot\dim_{\BF_p}(A'/A'a') 
\ =\ - \rank(\phi')\cdot\deg(\infty')\cdot\ord_{\infty'}(a')$$
for all $a'\in A'$. In other words $\deg_\tau(\phi'_{a'})=-v(a')$, where $v:=\rank(\phi')\cdot\deg(\infty')\cdot\ord_{\infty'}$ is a valuation on $F'$ equivalent to~$\ord_{\infty'}$. 
Now the slope of $\min_f$ with respect to $v$ is $-v(f)$; hence for each $i$ we have 
$$\deg_\tau(\phi_{a_i})\ =\ \deg_\tau(\phi'_{a_i})\ =\ -v(a_i)\ \le\ -i\cdot v(f)\ =\ i\cdot\deg_\tau(f),$$
with equality for $i=k$.
\end{Proof}

\begin{Prop}\label{MinPolEndEff}
For any $f\in\End_K(\phi)$, one can effectively compute $\min_f$.
% the minimal polynomial of $f$ over~$F$.
\end{Prop}

\begin{Proof}
Since $F[f]$ is a commutative subfield of $\End_K^{\circ}(\phi)$, the degree $k$ of $\min_f$ divides $\rank(\phi)$. Write $\min_f(X) = X^k+a_1X^{k-1}+\ldots+a_k$ with $a_i\in A$. By the definition of $\rank(\phi)$ and Proposition \ref{MinCoeff}, for each $1\le i\le k$ we then have $a_i=0$ or
%$$\rank(\phi)\cdot\dim_{\BF_p}(A/Aa_i) 
%= \deg_\tau(\phi_{a_i})\le i\cdot\deg_\tau(f)$$
$$\dim_{\BF_p}(A/Aa_i) 
= \frac{\deg_\tau(\phi_{a_i})}{\rank(\phi)}
\le \frac{i\cdot\deg_\tau(f)}{\rank(\phi)}.$$
% By Lemma \ref{AdLem}
Thus each $a_i$ lies in a finite set that can be effectively determined. For each choice of candidates for the $a_i$ one can effectively compute $f^k+\phi_{a_1}f^{k-1}+\cdots+\phi_{a_k}$ in $K[\tau]$ and decide whether it is zero. Thus $\min_f$ can be computed by letting $k$ run through the divisors of $r$ in ascending order and checking all possible choices of coefficients~$a_i$.
\end{Proof}

\begin{Prop}\label{KFinEff}
If $K$ is finite, one can effectively 
\begin{enumerate}
\item[(a)] find a finite separable extension $K'$ of $K$ such that $\End_{K'}(\phi)=\End_{K^\sep}(\phi)$, 
\item[(b)] compute the dimensions of $\End_{K^\sep}^\circ(\phi)$ and of its center over~$F$, and
\item[(c)] describe the center of $\End_{K^\sep}^\circ(\phi)$ as an abstract field extension of~$F$.
\end{enumerate}
\end{Prop}

\begin{Proof}
Set $\Frob_K := \tau^{[K/\BF_p]}$. Then for any finite extension $K_n/K$ of degree~$n$, the center of $\End_{K_n}^\circ(\phi)$ is a finite field extension of $F$ that is generated by $\Frob_{K_n} := \tau^{[K_n/\BF_p]} =\Frob_K^n$. As a special case of Proposition \ref{MinPolEndEff} one can effectively determine the minimal polynomial of $\Frob_K$ over~$F$. Using Lemma \ref{FieldExtLem}, one can therefore effectively find an integer $m\ge1$ such that $\bigcap_{n\ge1}F(\Frob_K^n) = F(\Frob_K^m)$. Whenever $K_m\subset K_n$, it follows that $\Frob_{K_m}=\Frob_K^m$ lies in the center of $\End_{K_n}^\circ(\phi)$. In particular, every element of $\End_{K_n}(\phi)$ commutes with $\Frob_{K_m}$ and is therefore defined over~$K_m$; in other words we have $\End_{K_n}(\phi)=\End_{K_m}(\phi)$. Varying $n$ we deduce that $\End_{K^\sep}(\phi) = \End_{K_m}(\phi)$, proving (a).

By Lemma \ref{FieldExtLem} one can also effectively calculate the minimal polynomial of $\Frob_{K_m}=\Frob_K^m$ over~$F$. 
This in turn determines the center $F(\Frob_{K_m})$ as an abstract field extension of~$F$, proving (c). In particular, its dimension $d:=\dim_F(F(\Frob_{K_m}))$ can be effectively computed. Moreover, let $e^2$ be the dimension of $\End_{K_m}^\circ(\phi)$ over $F(\Frob_{K_m})$. Since $K_m$ is finite, we then have $de=\rank(\phi)$. Thus $e$ can be effectively computed from $d$ and $\rank(\phi)$, which implies (b).
\end{Proof}

\begin{Prop}\label{ModelEff}
One can effectively construct a normal integral domain $R\subset K$ which is finitely generated over~$\BF_p$ with $\Quot(R)=K$, such that $\phi$ extends to a Drinfeld $A$-module over $\Spec R$. 
\end{Prop}

\begin{Proof}
By assumption $K$ is given as the fraction field of a finitely generated integral domain~$R$. For all elements $a$ of a finite set of non-zero generators of the $\BF_p$-algebra~$A$, adjoin to $R$ all coefficients of $\phi_a\in K[\tau]$ as well as the inverses of their highest coefficients. Then $R$ is still a finitely generated integral domain with $\Quot(R)=K$, and $\phi$ extends to a Drinfeld $A$-module over $\Spec R$. Finally, replace $R$ by its normalization. By our computer algebra prerequisites, all these operations can be carried out effectively.
\end{Proof}

\begin{Prop}\label{FrobEff}
For any maximal ideal $\Fm\subset R$ one can effectively 
\begin{enumerate}
\item[(a)] determine the  minimal and characteristic polynomials $\min_\Fm$ and $\charact_\Fm$ of $\Frob_\Fm$, and
% in $A[X]$, and
\item[(b)] decide whether $\End_{k_\Fm^\sep}^\circ(\phi_\Fm)$ is commutative and in that case describe it as an abstract field extension of~$F$.
\end{enumerate}
\end{Prop}

\begin{Proof}
For any $\Fm$ we have an explicit description of the reduction $\phi_\Fm$ over the finite residue field~$k_\Fm$. By Proposition \ref{MinPolEndEff} we can therefore effectively determine the associated minimal polynomial of $\Frob_\Fm := \tau^{[k_\Fm/\BF_p]}$ over~$F$. Consequently we can also effectively determine the characteristic polynomial $\charact_\Fm(X) := \min_\Fm(X)^{\rank(\phi)/\deg(\min_\Fm)}$, proving (a).
Part (b) is a direct consequence of Proposition \ref{KFinEff} (b) and (c).
\end{Proof}

%%%%%%%%%%%%%%%%%%%%%%%%%%%%%%%%%%%%%%%%%%%%%%%%%%%%%%%%%%%%%%%%%%%%%%%%%%%
%\newpage
\section{Searching for endomorphisms}
\label{SearchEndos}

We keep $(A,K,\phi)$ as before.

\begin{Prop}\label{EndodEff}
For any integer $d\ge0$, one can effectively determine a finite separable extension $K'$ of~$K$, such that all elements of $\End_{K^\sep}(\phi)$ of degree $d$ in~$\tau$ are defined over~$K'$, and determine these endomorphisms.
\end{Prop}

\begin{Proof}
Choose a finite set $S$ of generators of the $\BF_p$-algebra~$A$. Then an element $u\in K[\tau]$ is an endomorphism of $\phi$ if and only if $u\phi_a=\phi_au$ for all $a\in S$. With the Ansatz $u=\sum_{i=0}^{d}u_i\tau^i$ these equations amount to finitely many polynomial equations in the coefficients~$u_i$. We also know that there are at most finitely many solutions. By our computer algebra prerequisites, one can therefore effectively describe all these solutions and a common field of definition $K'$ for them. Since all endomorphisms over $\oK$ are already defined over~$K^\sep$, one can choose $K'$ separable over~$K$.
\end{Proof}

\medskip
Next consider the natural $A$-algebra homomorphism
\UseTheoremCounterForNextEquation
\begin{equation}\label{LowestCoeffMap}
\textstyle
\End_{K^\sep}(\phi)\to K^\sep, \;u=\sum_i u_i\tau^i\mapsto u_0.
\end{equation} 

\begin{Prop}\label{EndodCoeff}
Assume that $\phi$ has generic characteristic. Then:
\begin{enumerate}
\item[(a)] The homomorphism (\ref{LowestCoeffMap}) is injective.
\item[(b)] For any $u_0\in K^\sep$, one can effectively decide whether there exists $u\in \End_{K^{\sep}}(\phi)$ with lowest coefficient $u_0$ and, if so, find it.
\item[(c)] If $u\in \End_{K^{\sep}}(\phi)$ has lowest coefficient $u_0\in K$, then $u\in K[\tau]$.
\end{enumerate}
\end{Prop}

\begin{Proof}
By assumption the homomorphism is injective on~$A$. Since $\End_{K^\sep}(\phi)$ is an integral ring extension of~$A$, the homomorphism is  therefore injective. 

Suppose that $u_0\in K^\sep$ is the lowest coefficient of some element $u\in \End_{K^{\sep}}(\phi)$. Let $\min_u(X)\in A[X]$ be the minimal polynomial of $u$ over~$A$. Then $\min_u(u_0)=0$, hence $u_0$ is integral algebraic over~$A$, and by injectivity $\min_u$ is also the minimal polynomial of $u_0$ over~$A$. Write $\min_u(X) = X^k+a_1X^{k-1}+\ldots+a_k$ with $a_i\in A$. Then
$\deg_\tau(u) = \deg_\tau(\phi_{a_k})/k$ by Proposition \ref{MinCoeff}. Here the right hand side depends only on the joint minimal polynomial of $u$ and $u_0$; hence $\deg_\tau(u)$ is uniquely determined by~$u_0$. 

Now consider an arbitrary element $u_0\in K^\sep$. By our computer algebra prerequisites one can effectively decide whether $u_0$ is algebraic over $F$ and, if so, determine its minimal polynomial $\min_{u_0}$ over~$F$. In particular, one can decide whether $u_0$ is integral over~$A$, which is necessary for it to be the lowest coefficient of an endomorphism. If so, write $\min_{u_0}(X) = X^k+a_1X^{k-1}+\ldots+a_k$ with $a_i\in A$; then we can also compute the number $d := \deg_\tau(\phi_{a_k})/k$. By the above arguments any endomorphism $u\in \End_{K^{\sep}}(\phi)$ with lowest coefficient $u_0$ then satisfies $\deg_\tau(u)=d$. 

To see whether such $u$ actually exists we use the ansatz $u=\sum_{i=0}^{d}u_i\tau^i$ with $u_i\in K^\sep$. Pick an arbitrary non-constant $t\in A$ and write $\phi_t=\sum_{j=0}^nx_j\tau^j$. Expand the equation $u\phi_t=\phi_tu$ in the form
$$\sum_{i,j} u_ix_j^{p^i}\tau^{i+j}
\ =\ \sum_i u_i\tau^i\sum_jx_j\tau^j
\ =\ \sum_j x_j\tau^j\sum_iu_i\tau^i
\ =\ \sum_{i,j} x_ju_i^{p^j}\tau^{i+j}.$$ 
Comparing the coefficients of $\tau^\ell$ for each $0<\ell\le d$ yields an expression for ${(x_0-x_0^{p^\ell})u_\ell}$ as a polynomial in $x_1,\ldots,x_n$ and $u_0,\ldots,u_{\ell-1}$. Since $t$ is non-constant and $\phi$ has generic characteristic, the element $x_0$ is transcendental over~$\BF_p$. Thus $x_0-x_0^{p^\ell}\not=0$, and so we obtain explicit recursion relations for all $u_\ell$ in terms of the constant coefficient~$u_0$. Plugging these into the equations obtained by comparing the coefficients of $\tau^\ell$ for $\ell>d$ then yields explicit polynomial equations in $u_0$ over~$K$. These equations are fulfilled if and only if the desired $u$ exists, and in that case we can effectively find it. This proves (b).

Finally, the recursion relations show that $u\in K(u_0)[\tau]$ if it exists. In particular this implies (c).
\end{Proof}

\medskip
Now we fix a non-constant element $t\in A$ and set $\delta_t:=\deg_\tau(\phi_t)>0$. 
We view $K^\sep[\tau]$ as an $\BF_p[t]$-module via the multiplication $(a,u)\mapsto\phi_au$. Note that $K^\sep[\tau]$ is torsion free, so every finitely generated $\BF_p[t]$-submodule is free.

\begin{Def}\label{DefOrthBas}
We call a sequence of non-zero elements $m_1,\ldots,m_n$ of $K^\sep[\tau]$ \emph{orthogonal} if for all $a_1,\ldots,a_n\in\BF_p[t]$ we have 
$$\deg_\tau\Bigl(\,\sum_{i=1}^n\phi_{a_i}m_i\Bigr)\ =\ 
\max\bigl\{\deg_\tau(\phi_{a_i}m_i)\bigm|1\le i\le n\bigr\}.$$
\end{Def}

Here $\deg_\tau(\phi_{a_i}m_i)=\deg_t(a_i)\cdot\delta_t+\deg_\tau(m_i)$. In particular, we have $\deg_\tau(\phi_{a_i}m_i)\allowbreak=-\infty$ if and only if $\deg_t(a_i)=-\infty$ if and only if $a_i=0$. Thus the definition permits $\sum_{i=1}^n\phi_{a_i}m_i$ to be zero only if all $a_i$ are zero; hence any orthogonal sequence is $\BF_p[t]$-linearly independent.

\begin{Prop}\label{OrthMember}
Given a finite extension $K'$ of~$K$ and an $\BF_p[t]$-submodule $M$ of $K'[\tau]$ with an orthogonal basis $m_1,\ldots,m_n$, for any $f\in K'[\tau]$, one can effectively decide whether $f\in M$ and if so, one can effectively compute its coefficients with respect to that basis.  	
\end{Prop}

\begin{Proof}
If $f=\sum_i\phi_{a_i}m_i$ with $a_i\in\BF_p[t]$, by Definition \ref{DefOrthBas} we must have ${\deg_t(a_i)\cdot\delta_t}\allowbreak+\allowbreak\deg_\tau(m_i)\le \deg_\tau(f)$ for all~$i$.
Thus each $\deg_t(a_i)$ is bounded by an explicit number. Writing each $a_i=\sum_ja_{ij}t^j$ with finitely many $a_{ij}\in\BF_p$ to be determined, the equation $f=\sum_i\phi_{a_i}m_i$ is equivalent to finitely many linear equations in the $a_{ij}$ with coefficients in~$K^\sep$. By our computer algebra prerequisites, one can effectively decide whether these have a solution in $\BF_p$ and if so, find it.
\end{Proof}

%\medskip
%Identifying $A$ with its image under the embedding $\phi\colon A\into K[\tau]$, we can apply the concept to elements of~$A$:
%
%\begin{Prop}\label{AOrthBasis}
%One can effectively find an orthogonal basis of $A$ over $\BF_p[t]$.
%\end{Prop}
%
%\begin{Proof}
%By the definition of $\rank(\phi)$, for all $a\in A$ we have 
%%$\deg_\tau(\phi_a)=\rank(\phi)\cdot\dim_{\BF_p}(A/Aa)$. 
%$$\deg_\tau(\phi_a) 
%\ =\ \rank(\phi)\cdot\dim_{\BF_p}(A/Aa) 
%\ =\ - \rank(\phi)\cdot\deg(\infty)\cdot\ord_\infty(a),$$
%where $\ord_\infty$ denotes the normalized valuation on $F$ associated to~$\infty$ and $\deg(\infty)$ the degree of its residue field over~$\BF_p$. Thus a basis $m_1,\ldots,m_n$ of $A$ over $\BF_p[t]$ is orthogonal if and only if for all $a_1,\ldots,a_n\in\BF_p[t]$ we have 
%$$\ord_\infty\Bigl(\,\sum_{i=1}^n a_im_i\Bigr)\ =\ 
%\min\bigl\{\ord_\infty(a_im_i)\bigm|1\le i\le n\bigr\}.$$
%... ? ? ? ...
%\end{Proof}

\medskip
For any integer $d$ let $M_d$ denote the $\BF_p[t]$-submodule of $\End_{K^\sep}(\phi)\subset K^\sep[\tau]$ generated by all $u \in \End_{K^\sep}(\phi)$ with $\deg_\tau(u)\le d$. As there are only finitely many generators, this module is finitely generated. Recall that all finite separable extensions of $K$ are tacitly assumed to be contained in~$K^\sep$.

\begin{Prop}\label{OrthFindStep}
For any $d\ge0$, any finite separable extension $K'$ of~$K$, and any elements $m_1,\ldots,m_n \in \End_{K'}(\phi)$ which form an orthogonal basis of $M_{d-1}$, one can effectively find $n'\ge n$, a finite separable extension $K''$ of~$K'$, and elements $m_{n+1},\ldots,m_{n'}\in \End_{K''}(\phi)$ of degree $d$ such that $m_1,\ldots,m_{n'}$ is an orthogonal basis of~$M_d$.
\end{Prop}

\begin{Proof}
First observe that, since $M_{d-1}$ is generated by elements of degree $\le d-1$, Definition \ref{DefOrthBas} implies that $\deg_\tau(m_i)\le d-1$ for all $1\le i\le n$. 

Next, applying Proposition \ref{EndodEff} over~$K'$, we can effectively find a finite separable extension $K''$ of~$K'$ and elements $f_1,\ldots,f_k\in\End_{K''}(\varphi)$ which make up all elements of $\End_{K^\sep}(\varphi)$ of degree $d$ in~$\tau$. 
By definition, these together with $m_1,\ldots,m_n$ generate~$M_d$.
Using Proposition \ref{OrthMember}, for all $(\alpha_1,\ldots,\alpha_k)\in\BF_p^k\setminus\{(0,\ldots,0)\}$ we can check whether $\sum_{j=1}^k\alpha_jf_j$ already lies in $M_{d-1}$. 
If it does, we can remove one generator $f_j$ without changing the module~$M_d$. After finitely many such operations and reordering $f_1,\ldots,f_k$, we may assume that we have constructed an integer $0\le\ell\le k$ such that $M_d$ is already generated by $m_1,\ldots,m_n,f_1,\ldots,f_\ell$ and that:
\UseTheoremCounterForNextEquation
\begin{equation} \label{OrthFindRed}
\hbox{For all $(\alpha_1,\ldots,\alpha_\ell)\in\BF_p^\ell\setminus\{(0,\ldots,0)\}$ we have $\sum_{j=1}^\ell\alpha_jf_j\not\in M_{d-1}$.}
\end{equation}
%for all $(\alpha_1,\ldots,\alpha_\ell)\in\BF_p^\ell\setminus\{(0,\ldots,0)\}$ we have $\sum_{j=1}^\ell\alpha_jf_j\not\in M_{d-1}$.
We claim that then $m_1,\ldots,m_n,f_1,\ldots,f_\ell$ are orthogonal.

To see this, set $n':=n+\ell$ and $m_{n+j}:=f_j$ for all $1\le j\le\ell$.
For the sake of contradiction consider $a_1,\ldots,a_{n'}\in\BF_p[t]$ such that
\UseTheoremCounterForNextEquation
\begin{equation} \label{CancelDeg}
\deg_\tau\Bigl(\,\sum_{i=1}^{n'}\phi_{a_i}m_i\Bigr)
\ <\ D\ :=\ 
\max\bigl\{\deg_\tau(\phi_{a_i}m_i)\bigm|1\le i\le n'\bigr\}.
\end{equation}
Then the maximum $D$ is attained for some $n<i\le n'$, because otherwise the strict inequality would also hold with $a_{n+1},\ldots,a_{n'}$ replaced by~$0$, contradicting the orthogonality of $m_1,\ldots,m_n$. Thus $D\ge d$.

For all $1\le i\le n'$ we therefore have $D\ge\deg_\tau(m_i)$. Write $D=s_i\delta_t+\deg_\tau(m_i)$, so that $\deg_t(a_i)\le s_i$. Dropping from $a_i\in\BF_p[t]$ all terms of degree $<s_i$ in~$t$ then does not change the inequality (\ref{CancelDeg}). After doing this for all $i$ we may assume that each $a_i=\alpha_i t^{s_i}$ for some $\alpha_i\in\BF_p$, with $\alpha_i=0$ if $s_i\not\in\BZ$. 

If $D>d$, we have $s_i>0$ for all~$i$. This now means that all $a_i$ are divisible by~$t$. By the additivity of $\deg_\tau$ in products, the inequality (\ref{CancelDeg}) still holds after replacing each $a_i$ by $t^{-1}a_i$. After repeating this a finite number of times, we can assume that $D=d$. 

Then $s_i=0$ for all $n<i\le n'$. For these we then have $a_i=:\alpha_i\in\BF_p$. Also, since the maximum in (\ref{CancelDeg}) is attained at some $n<i\le n'$, at least one of these $\alpha_i$ is non-zero.
On the other hand the inequality (\ref{CancelDeg}) now means that $\deg_\tau(m)<d$ for $m:=\sum_{i=1}^{n'}\phi_{a_i}m_i$. Since $m\in\End_{K^\sep}(\phi)$, it follows that $m\in M_{d-1}$. This in turn implies that 
$$\sum_{j=1}^\ell\alpha_{n+j}f_j
\ =\ \sum_{i=n+1}^{n'}\phi_{a_i}m_i
\ =\ m-\sum_{i=1}^n\phi_{a_i}m_i
\ \in\ M_{d-1}.$$
But this is a contradiction to (\ref{OrthFindRed}).

This proves that $m_1,\ldots,m_n,f_1,\ldots,f_\ell$ are orthogonal. 
As any orthogonal sequence is $\BF_p[t]$-linearly independent, and these elements generate $M_d$ as an $\BF_p[t]$-module; they form an orthogonal basis of $M_d$, as desired.
\end{Proof}

\begin{Prop}\label{OrthFind}
For any $d$ one can effectively determine a finite separable extension $K'$ of $K$ and elements $m_1,\ldots,m_n \in \End_{K'}(\phi)$ which form an orthogonal basis of~$M_d$.
\end{Prop}

\begin{Proof}
For $d<0$ we have $M_d=0$ and the assertion is trivial. By induction on $d$ the assertion thus follows from Proposition \ref{OrthFindStep}.
\end{Proof}

\medskip
Since $\End_{K^\sep}(\phi)$ is finitely generated as a module over $A$ and hence also over $\BF_p[t]$, we have $\End_{K^\sep}(\phi)=M_d$ for all $d\gg0$. Letting the procedure in Proposition \ref{OrthFindStep} run inductively for $d\to\infty$, we will therefore eventually find an orthogonal basis of $\End_{K^\sep}(\phi)$ over $\BF_p[t]$. However, knowing when we have reached that stage requires additional information. By the following proposition it suffices to know when $M_d$ has the same rank as $\End_{K^\sep}(\phi)$:

\begin{Prop}\label{OrthRankEnough}
For any~$d$, if $M_d$ has finite index in $\End_{K^\sep}(\phi)$, it is equal to $\End_{K^\sep}(\phi)$.
\end{Prop}

\begin{Proof}
Let $m_1,\ldots,m_n$ be an orthogonal basis of $\End_{K^\sep}(\phi)$. Since $M_d$ is generated by elements of degree $\le d$ and has finite index in $\End_{K^\sep}(\phi)$, for any $1\le i\le n$ there exists an element $m\in M_d$ of degree $\le d$, in whose expansion $m=\sum_{j=1}^n\phi_{a_j}m_j$ with $a_j\in\BF_p[t]$ the coefficient $a_i$ is non-zero. By orthogonality we then have 
$$d \ge\deg_\tau(m) \ge \deg_\tau(\phi_{a_i}m_i) \ge \deg_\tau(m_i).$$
Thus the generators $m_i$ of $\End_{K^\sep}(\phi)$ already lie in~$M_d$; hence $\End_{K^\sep}(\phi)=M_d$.
\end{Proof}

\begin{Prop}\label{RankEff}
If the rank of $\End_{K^\sep}(\phi)$ over $\BF_p[t]$ is given, one can effectively determine a finite separable extension $K'$ of $K$ and elements $m_1,\ldots,m_n \in \End_{K'}(\phi)$ which form an orthogonal basis of $\End_{K^\sep}(\phi)$.
\end{Prop}

\begin{Proof}
Let the procedure in Proposition \ref{OrthFindStep} run inductively for $d\to\infty$ until $\rank(M_d)=\rank(\End_{K^\sep}(\phi'))$. Then we are finished by Proposition \ref{OrthRankEnough}.
\end{Proof}

\begin{Rem} \label{OrthRem}
\rm Given a finite extension $K'$ of~$K$ and an orthogonal basis $m_1,\ldots,m_n$ of $\End_{K'}(\phi)$, we can effectively answer various elementary questions about this endomorphism ring. Namely, using Proposition \ref{OrthMember} we can write each product $m_im_j$ as an $\BF_p[t]$-linear combination of $m_1,\ldots,m_n$ and thus describe $\End_{K'}(\phi)$ via a multiplication table. By solving linear equations over $\BF_p[t]$ we can then, for instance, explicitly determine the commutant of any element, or of any finite number of elements, of $\End_{K'}(\phi)$. In particular, we can explicitly determine the center of $\End_{K'}(\phi)$. 
We can also find a finite presentation of $\End_{K'}(\phi)$ as an $A$-module or as an $A$-algebra.
Thus we can say that \emph{we know} $\End_{K'}(\phi)$. 
\end{Rem}

\begin{Rem}\label{IsogEndEffRem}
\rm Suppose that we are given a Drinfeld $A$-module $\phi'$ isogenous to $\phi$ and an orthogonal basis of $\End_{K^\sep}(\phi')$. 
Then $\End_{K^\sep}(\phi)$ and $\End_{K^\sep}(\phi')$ have the same rank over $\BF_p[t]$; hence we can effectively find an orthogonal basis of $\End_{K^\sep}(\phi)$ using Proposition \ref{RankEff}. But this seems wasteful. It is probably more economical to use given isogenies $\phi\to\phi'\to\phi$ 
%explicitly to find a bound for the degree $d$ or 
to find the endomorphisms of $\phi$ explicitly from those of~$\phi'$.
\end{Rem}

\begin{Prop}\label{AllEndosKFinite}
If $K$ is finite, for any non-constant element $t\in A$ one can effectively find an orthogonal basis of $\End_{K}(\phi)$ over~$\BF_p[t]$.
\end{Prop}

\begin{Proof}
Choose a basis $\{x_i\mid 1\le i\le n\}$ of $K$ over~$\BF_p$. Then ${\{x_i\tau^j\mid 1\le i\le n,\; 0\le j<n\}}$ is a basis of $K[\tau]$ as a free module over its center $\BF_p[\tau^n]$. For any $a\in A$ one can explicitly write $\phi_a$ as an $\BF_p[\tau^n]$-linear combination of this basis. Doing this for finitely many generators of $A$ as an $\BF_p$-algebra, one can then explicitly determine their joint commutant in $K[\tau]$ by linear algebra over $\BF_p[\tau^n]$. By definition this commutant is precisely $\End_K(\phi)$. Thus one has an explicit basis of $\End_K(\phi)$ as a module over $\BF_p[\tau^n]$. 

By another explicit calculation over $\BF_p[\tau^n]$ one can determine the action of $\phi_t$ on this basis; in other words, one can determine $\End_K(\phi)$ as a module over the polynomial ring in two variables $\BF_p[\tau^n,t]$. That in turn yields a basis of $\End_K(\phi)$ as a module over $\BF_p[t]$. If this calculation is done carefully, say with a suitable Gr\"obner basis, the basis is already orthogonal; otherwise an orthogonal basis can be obtained from this one as in the proof of Proposition \ref{OrthFindStep}.
\end{Proof}

%%%%%%%%%%%%%%%%%%%%%%%%%%%%%%%%%%%%%%%%%%%%%%%%%%%%%%%%%%%%%%%%%%%%%%%%%%%
%\newpage
\section{Main algorithms}
\label{MainAlg}

As before we fix a Drinfeld $A$-module $\phi\colon A\to K[\tau]$ over a finitely generated field~$K$. In this section we show that one can effectively determine the endomorphism ring $\End_{K^\sep}(\phi)$ and, in the non-isotrivial special characteristic case, the admissible coefficient rings $B\subset A'$ and the Drinfeld $A'$-module $\phi'$ from Theorem \ref{EndGalSpec1} and the endomorphism ring $\End_{K^\sep}(\phi'|B)$.

\medskip
We begin by determining whether $\phi$ has generic or special characteristic and whether it is isotrivial, using Proposition \ref{CharIsoEff}.
If $\phi$ is not isotrivial, we first find a maximal commutative subring of $\End_{K^\sep}(\phi)$, whose normalization will be the ring $A'$ below. We say that we \emph{go~up} with the coefficient ring.

\begin{Prop}\label{GoingUp}
If $\phi$ is not isotrivial, one can effectively find a finite separable extension $K'$ of~$K$, an admissible coefficient ring $A'$ containing~$A$, a Drinfeld $A'$-module $\phi'$ over~$K'$, and an isogeny $h\colon\phi\to\phi'|A$ over~$K'$, such that $\End_{K^\sep}(\phi')=A'$.
%One can effectively find a finite separable extension $K'$ of~$K$, an admissible coefficient ring $A'$ containing~$A$, a Drinfeld $A'$-module $\phi'$ over~$K'$, isogenies $f\colon\phi\to\phi'|A$ and $g\colon\phi'|A\to\phi$ over~$K'$, and an element $a\in A\setminus\{0\}$ with $gf=\phi_a$, such that $\End_{K^{\prime\sep}}(\phi')=A'$.
\end{Prop}

\begin{Proof}
Set $(A_0,K_0,\phi_0) := (A,K,\phi)$ and $n:=0$, and start processes (a) and (b) in parallel.

\medskip{\it Process (a): Find endomorphisms:} 
For each $d\ge0$ use Proposition \ref{EndodEff} to find all endomorphisms $f\in\End_{K^\sep}(\phi_m)$ of degree~$d$. For any such $f$ use Proposition \ref{MinPolEndEff} to check whether $f$ is scalar or not. 
If a non-scalar $f$ is found, choose a finite separable extension $K_{n+1}$ of $K_n$ over which $f$ is defined. Let $A_{n+1}$ be the normalization of $A_n[f]$, and using Proposition \ref{IsogEff} choose a Drinfeld $A_{n+1}$-module $\phi_{n+1}$ and an isogeny $\phi_n\to\phi_{n+1}|A_n$ over $K_{n+1}$. Then kill process (b), set $n:=n+1$, and restart both processes (a) and (b).

\medskip{\it Process (b): Find Frobeniuses:} 
Using Proposition \ref{ModelEff} choose a finitely generated normal integral domain $R_n\subset K_n$ with $\Quot(R_n)=K_n$ over which $\phi_n$ has good reduction. For each maximal ideal $\Fm\subset R_n$ use Proposition \ref{FrobEff} (b) to decide whether $F_{n,\Fm} := \End_{k_\Fm^\sep}^\circ(\phi_{n,\Fm})$ is commutative, and in that case describe it as an abstract field extension of $F_n := \Quot(A_n)$.
As soon as a new such $F_{n,\Fm}$ is found, check whether it and some previously found $F_{n,\Fm'}$ are linearly disjoint over~$F_n$. 
If no, continue with the next~$\Fm$.
If yes, we know that $\End_{K^\sep}(\phi_n)=A_n$. 
Then kill process (a), set $(A',K',\phi'):=(A_n,K_n,\phi_n)$, combine all isogenies from process (a) to an isogeny $h\colon\phi\to\phi'|A$ over~$K'$, and stop.

\medskip{\it Effectivity:} 
Process (a) constructs a sequence of Drinfeld modules of strictly decreasing rank. Thus eventually it continues forever with the same $(A_n,K_n,\phi_n)$. In that case we have $\End_{K^\sep}(\phi_n)=A_n$. 
Process (b) cannot terminate before that, because the rings $F_{n,\Fm}$ all contain a subring isomorphic to $\End_{K^\sep}^\circ(\phi_n)$ over~$F_n$. But once $(A_n,K_n,\phi_n)$ remains constant, by Proposition \ref{X1} there exist maximal ideals $\Fm$ and $\Fn$ of $R_n$ such that $F_{n,\Fm}$ and $F_{n,\Fn}$ are linearly disjoint over~$F_n$. Thus process (b) terminates with a correct answer.
\end{Proof}

\begin{Var}\label{GoingUpVar}
\rm In process (b) of Proposition \ref{GoingUp}, instead of checking for linear disjointness, for each $\Fm$ make a list $\CL_\Fm$ of the finitely many isomorphism classes of field extensions $E/F_n$ with $\Hom_{F_n}(E,F_{n,\Fm})\not=\emptyset$. For any new $\Fm$ compare this list with all previously found lists $\CL_{\Fm'}$. If the intersection of these is the singleton $\{F_n\}$, kill process (a) and finish as before.
\end{Var}

\medskip
If $\phi$ is non-isotrivial of special characteristic, we must \emph{go~down} with the coefficient ring in order to discover more endomorphisms.

\begin{Prop}\label{GoingDown}
If $\phi$ is non-isotrivial of special characteristic, let $(A',K',\phi')$ be the data from Proposition \ref{GoingUp}. Then one can effectively find the admissible coefficient ring $B\subset A'$ from Theorem \ref{EndGalSpec1}, a non-constant element $t\in B$, a finite separable extension $K''$ of $K'$, and elements of $\End_{K''}(\phi'|B)$ which form  an orthogonal basis of $\End_{K^\sep}(\phi'|B)$ over~$\BF_p[t]$.
\end{Prop}

\begin{Proof}
Since $\phi'$ is non-isotrivial we have $\rank(\phi')\ge2$. 
If $\rank(\phi')=2$, then $B=\End_{K^\sep}(\phi'|B)=A'$ by Proposition \ref{Rank2}. 
In particular, for any non-constant element $t\in B$ we know the rank of $B$ over~$\BF_p[t]$; hence we can effectively find an orthogonal basis over $\BF_p[t]$ using Proposition \ref{RankEff}.
%NOTE: In practice one would much rather find an orthogonal basis of $B$ by calculation within $B$ alone, using Proposition \ref{AOrthBasis}, but we did not want to write out the extra arguments.

So assume that $\rank(\phi')>2$. Then $B$ is completely characterized by traces by Theorem \ref{EndGalSpec2}. Start process (a).

\medskip{\it Process (a): Find traces of Frobenius:} 
Using Proposition \ref{ModelEff} choose a finitely generated normal integral domain $R'\subset K'$ with $\Quot(R')=K'$ over which $\phi'$ has good reduction. Set $F' := \Quot(A')$ and $k:=0$ and $B_0 := \BF_p$. 
For each maximal ideal $\Fm'\subset R'$ use Proposition \ref{FrobEff} (a) to compute the characteristic polynomial $\charact_{\Fm'}=\sum_{i=0}^{r'}a_iX^i\in F'[X]$ of $\Frob_{\Fm'}$ associated to~$\phi'$. Using this, calculate the value $t_{\Fm'} = a_1a_{r'-1}/a_0 \in F'$ from (\ref{tmDef}).
If $t_{\Fm'}\not\in\Quot(B_k)$, determine $B_{k+1}:= A'\cap\Quot(B_k[t_{\Fm'}]) \subset B$ and set $k:=k+1$. 

Keep repeating this forever with all~$\Fm'$. 
The first time that $B_k$ becomes infinite, fix a non-constant element $t\in B_k$, set $B':=\BF_p[t]$, and start process (b) in parallel.

\medskip{\it Process (b): Find endomorphisms:}
For all integers $d$ let $M'_d$ denote the $B'$-submodule of $\End_{K^\sep}(\phi'|B')$ generated by all $u \in \End_{K^\sep}(\phi'|B')$ with $\deg_\tau(u)\le d$. Thus ${M'_{-1}=0}$ with the trivial orthogonal basis.
Using Proposition \ref{OrthFindStep} inductively, for every $d\ge0$ we can effectively construct a finite separable extension $K'_d$ of $K'$ and an orthogonal basis of~$M'_d$ contained in $\End_{K'_d}(\phi'|B')$.
If, with the current $B_k$ from process (a), we have
\UseTheoremCounterForNextEquation
\begin{equation}\label{GoingDownEq1}
\rank_{B'}(B_k)\cdot\rank_{B'}(M'_d)
\ =\ \rank_{B'}(A')^2,
\end{equation}
kill process (a) and stop, returning $B_k$ and $K'':=K'_d$ and the given orthogonal basis of~$M'_d$.

\medskip{\it Effectivity:} 
Process (a) produces an increasing sequence of normal subrings $B_k$ of~$A'$. By Theorem \ref{EndGalSpec2} this sequence eventually becomes stationary with $B_k=B$. In particular, from some point on $B_k$ is infinite and hence an admissible coefficient ring. Process (b) then produces an increasing sequence of $B'$-submodules $M'_d$ of $\End_{K^\sep}(\phi'|B')$ which eventually becomes stationary with $M'_d=\End_{K^\sep}(\phi'|B')$. By Theorem \ref{EndGalSpec1} (b) we have 
$\End_{K^\sep}(\phi'|B') \subset \End_{K^\sep}(\phi'|B)$ and hence 
$\End_{K^\sep}(\phi'|B') = \End_{K^\sep}(\phi'|B)$.
Thus at every step in process (b) we have
$$\rank_{B'}(B_k)\cdot\rank_{B'}(M'_d)
\ \le\ \rank_{B'}(B)\cdot\rank_{B'}(\End_{K^\sep}(\phi'|B)),$$
with equality for all $k,d\gg0$. But since $B$ is the center of $\End_{K^\sep}(\phi'|B)$ by Theorem \ref{EndGalSpec1} (a), and $A'$ is a maximal commutative subalgebra of $\End_{K^\sep}(\phi'|B)$, the right hand side of this inequality is equal to
$$\rank_{B'}(B)^2\cdot\rank_{B}(\End_{K^\sep}(\phi'|B))
\ =\ \rank_{B'}(B)^2\cdot\rank_{B}(A')^2
\ =\ \rank_{B'}(A')^2.$$
Thus at every step in process (b) we have
\UseTheoremCounterForNextEquation
\begin{equation}\label{GoingDownEq2}
\rank_{B'}(B_k)\cdot\rank_{B'}(M'_d)
\ \le\ \rank_{B'}(A')^2,
\end{equation}
with equality for all $k,d\gg0$. 
Comparing (\ref{GoingDownEq2}) with (\ref{GoingDownEq1}) shows that the process terminates and that upon termination we have $\rank_{B'}(B_k) = \rank_{B'}(B)$ and $\rank_{B'}(M'_d) = \rank_{B'}(\End_{K^\sep}(\phi'|B))$. The first of these equalities implies that $B_k=B$ because $B_k$ is normal, and the second implies that $M'_d = \End_{K^\sep}(\phi'|B)$ by Proposition \ref{OrthRankEnough}.
\end{Proof}

\begin{Var}\label{GoingDownVar}
\rm In process (a) of Proposition \ref{GoingDown}, in addition to $t_{\Fm'}$ adjoin all coefficients of the characteristic polynomial of $\Frob_{\Fm'}$ in the adjoint representation on $\End_{A'_{\Fp'}}(T_{\Fp'}(\phi'))$.
%for any maximal ideal $\Fp'$ different from the characteristic of~$\phi'_{\Fm'}$.
Like $t_{\Fm'}$ these coefficients can be computed directly from the characteristic polynomial of $\Frob_{\Fm'}$ on $T_{\Fp'}(\phi')$, and by the last sentence of Pink \cite[Thm.\;1.3]{PinkDrinSpec2}, they also lie in $E^\trad=E$. In this way one can probably generate $E$ faster.
\end{Var}

\begin{Thm}\label{RankKSep}
One can effectively determine the rank of $\End_{K^\sep}(\phi)$ over~$A$ and a finite separable extension $K''$ of $K$ with $\End_{K''}(\phi) = \End_{K^\sep}(\phi)$.
\end{Thm}

\begin{Proof}
First assume that $\phi$ has generic characteristic. Let $(A',K',\phi')$ be the data from Proposition \ref{GoingUp}. Then $\End_{K^\sep}(\phi)$ is commutative, hence so is $\End_{K^\sep}(\phi'|A)$; hence the latter is equal to $\End_{K^\sep}(\phi')=A'$. Thus the rank of $\End_{K^\sep}(\phi)$ over~$A$ is equal to that of $\End_{K^\sep}(\phi'|A)=A'$ and can be determined from the knowledge of~$A'$.
Also, since $\phi'$ and the isogeny $\phi\to\phi'|A$ are defined over~$K'$, it follows that $\End_{K^\sep}(\phi) = \End_{K'}(\phi)$. Thus the theorem holds with $K''=K'$. 

Next assume that $\phi$ is non-isotrivial of special characteristic. Let $(A',B,t,K'',\phi')$ be the data from Propositions \ref{GoingUp} and \ref{GoingDown} and abbreviate $S := \End_{K^\sep}(\phi'|A)$. Then by Theorem \ref{EndGalSpec1} (b) we have $S \subset \End_{K^\sep}(\phi'|B)$; hence $S$ is simply the commutant of $\phi'(A)$ in $\End_{K^\sep}(\phi'|B)$. Choose finitely many generators of $A$ as an $\BF_p$-algebra. By Proposition \ref{OrthMember} we can effectively express them as $\BF_p[t]$-linear combinations of the orthogonal basis from Proposition \ref{GoingDown}. By Remark \ref{OrthRem} we can therefore effectively determine $S$ as an $\BF_p[t]$-module. In particular, we can compute its rank over $\BF_p[t]$. Although $t$ does not necessarily lie in~$A$, we nevertheless have $A\cup\BF_p[t]\subset A'\subset S$. The multiplicativity of ranks thus implies that
$$\rank_A(S)\ =\ 
\rank_A(A')\cdot\rank_{A'}(S)\ =\ 
\rank_A(A')\cdot\frac{\rank_{\BF_p[t]}(S)}{\rank_{\BF_p[t]}(A')}.$$
Thus we can also compute the rank of $S$ over~$A$. As the rank of the endomorphism ring is invariant under isogenies, we have thereby computed the rank of $\End_{K^\sep}(\phi)$ over~$A$. 
Moreover, since $\phi'$ and the isogeny $\phi\to\phi'|A$ are defined over~$K''$, and $\End_{K^\sep}(\phi'|B) = \End_{K''}(\phi|B)$ by Proposition \ref{GoingDown}, it follows that $\End_{K^\sep}(\phi) = \End_{K''}(\phi)$, as desired.

Finally assume that $\phi$ is isotrivial. Pick a non-constant $t\in A$ and write $\phi_t=\sum_{i=0}^{n}x_i\tau^i$ with $x_n\not=0$. Choose a finite extension $K'$ of $K$ containing an element $y$ with $y^{1-q^n}=x_n$. Then $y^{-1}\phi_ty$ has the highest term~$\tau^n$. Since $\phi$ is isotrivial, as in the proof of Proposition \ref{CharIsoEff} it follows that the Drinfeld $A$-module $\psi := y^{-1}\phi y$ is now defined over a finite subfield $k$ of~$K'$. In fact, such $k$ can be described explicitly as the subfield generated by all coefficients of $\psi_a$ for a finite set of generators $a$ of $A$ as an $\BF_p$-algebra.
Applying Proposition \ref{KFinEff} to $(k,\psi)$ in place of $(K,\phi)$, one can then effectively find a finite separable extension $k'$ of~$k$ such that $\End_{k'}(\psi)=\End_{k^\sep}(\psi)$ and compute the rank of $\End_{k^\sep}(\psi)$ over~$A$. For any finite extension $K''$ of $K'$ containing a subfield isomorphic to $k'$ we then have $\End_{K''}(\phi) = \End_{K^\sep}(\phi)$ and $\rank_A(\End_{K^\sep}(\phi))=\rank_A(\End_{k^\sep}(\psi))$ and are done.
\end{Proof}

\begin{Thm}\label{AllEndosKSep}
For any non-constant element $t\in A$, one can effectively find a finite separable extension $K''$ of $K$ and elements of $\End_{K''}(\phi)$ which form an orthogonal basis of $\End_{K^\sep}(\phi)$ over~$\BF_p[t]$.
\end{Thm}

\begin{Proof}
{}From the knowledge of $A$ we can determine the rank of $A$ over~$\BF_p[t]$. Using Theorem \ref{RankKSep} we can therefore effectively determine the rank of $\End_{K^\sep}(\phi)$ over~$\BF_p[t]$. We can then find the desired data by Proposition \ref{RankEff}. (But in practice it might be more efficient to use the endomorphisms already found in Propositions \ref{GoingUp} and \ref{GoingDown}; compare Remark \ref{IsogEndEffRem}.)
\end{Proof}

\begin{Thm}\label{AllEndosK}
\begin{enumerate}
\item[(a)] One can effectively determine the rank of $\End_K(\phi)$ over~$A$.
\item[(b)] For any non-constant element $t\in A$, one can effectively find an orthogonal basis of $\End_{K}(\phi)$ over~$\BF_p[t]$.
\end{enumerate}
\end{Thm}

\begin{Proof}
Maybe this can be achieved by carrying out the whole program over~$K$ instead of~$K^\sep$, but we deduce it from the orthogonal basis in Theorem \ref{AllEndosKSep}, as follows. Let $m_1,\ldots,m_n\in\End_{K''}(\phi)$ be that basis. After replacing $K''$ by the subfield generated over $K$ by the coefficients of all~$m_i$, we can assume that $K''$ is galois over~$K$. For any $\sigma\in\Gal(K''/K)$ and any $i$ we can then express $\sigma(m_i)$ as an $\BF_p[t]$-linear combination of $m_1,\ldots,m_n$, using Proposition \ref{OrthMember}. In this way we can explicitly describe the action of $\Gal(K''/K)$ on $\End_{K''}(\phi)$. By solving linear equations over $\BF_p[t]$, we can then compute a basis of the submodule of invariants, which is precisely $\End_K(\phi)$. In particular, we can determine the rank of $\End_K(\phi)$ over $\BF_p[t]$, and hence also over~$A$. With a little more care we can make the basis orthogonal: Using the given orthogonal basis of $\End_{K''}(\phi)$, for any $d$ we can effectively find all elements of 
$\End_K(\phi)$ of degree $d$, and can then find an orthogonal basis as in the proof of Proposition \ref{OrthFindStep}.
\end{Proof}

\begin{Prop}\label{EndosOnTate}
One can effectively determine the isomorphism class of $T_\ad(\phi)$ as a module over $\End_{K^\sep}(\phi)\otimes_AA_\ad$.
\end{Prop}

\begin{Proof}
Recall that $S := \End_{K^\sep}^\circ(\phi)$ is a finite dimensional division algebra over~$F$. Let $Z$ denote its center. 
For any maximal ideal $\Fp\not=\Fp_0$ of~$A$ we then have $Z_\Fp := Z\otimes_FF_\Fp \cong \prod_\FP Z_\FP$, where the product is extended over all primes $\FP$ of $Z$ above~$\Fp$. Also $S_\Fp := S\otimes_FF_\Fp\cong \prod_\FP S_\FP$ where each $S_\FP:= S\otimes_ZZ_\FP$ is a central simple algebra over the field~$Z_\FP$. The isomorphism class of any $S_\FP$-module is therefore determined by its dimension over~$Z_\FP$. 

The rational $\Fp$-adic Tate module of $\phi$ is the $F_\Fp$-vector space $V_\Fp(\phi) := T_\Fp(\phi)\otimes_{A_\Fp}F_\Fp$. Under the above decomposition of $Z_\Fp$ it has the natural decomposition
$$V_\Fp(\phi) \ \cong\ V_\Fp(\phi'|A) 
\ =\ \prod_\FP V_\FP(\phi'|A'\cap Z).$$
Here each factor $V_\FP(\phi'|A'\cap Z)$ is a $Z_\FP$-vector space of dimension the rank of $\phi'|A'\cap Z$, which is $\rank(\phi)/[Z/F]$. 
Thus $V_\Fp(\phi)$ is a free module over $Z_\Fp$ of rank $\rank(\phi)/[Z/F]$. This therefore determines the isomorphism class of $V_\Fp(\phi)$ as a module over~$S_\Fp$.

Next, Theorem \ref{AllEndosKSep} and Remark \ref{OrthRem} yield a finite separable extension $K''$ of $K$ such that $\End_{K''}(\phi)=\End_{K^\sep}(\phi)$ and an explicit presentation of $\End_{K''}(\phi)$ as an $A$-algebra. Using the reduced trace of~$S$, one can find a maximal $A$-order $M\subset S$ containing $\End_{K^\sep}(\phi)$. 
For any maximal ideal $\Fp\not=\Fp_0$ of~$A$ the ring $M_\Fp := M\otimes_AA_\Fp$ is then a maximal order in~$S_\Fp$. In fact, we have $M_\Fp \cong \prod_\FP M_\FP$ for maximal orders $M_\FP$ in~$S_\FP$. Each $M_\FP$ is a matrix ring over a maximal $A_\Fp$-order $M_\FP'$ in a division algebra over~$F_\Fp$. Here $M_\FP'$ is a (possibly non-commutative) discrete valuation ring, and so any finitely generated torsion free $M_\FP'$-module is free. It follows that any finitely generated $A_\Fp$-torsion free $M_\Fp$-module is projective, and so its isomorphism class is determined by the $S_\Fp$-module obtained by base extension.
In particular, we therefore know the isomorphism class of the $M_\Fp$-submodule $\tTp$ of $V_\Fp(\phi)$ that is generated by $T_\Fp(\phi)$.

By analyzing the finite $A$-module $M/\End_{K^\sep}(\phi)$ one can construct a non-zero element $a\in A$ such that $a\cdot M \subset\End_{K^\sep}(\phi)$. 
For any maximal ideal $\Fp\not=\Fp_0$ with $\Fp\nmid a$ we then have $\End_{K^\sep}(\phi)\otimes_AA_\Fp = M_\Fp$ and hence $\tTp=T_\Fp(\phi)$, which determines the isomorphism class of $T_\Fp(\phi)$ as a module over $\End_{K^\sep}(\phi)\otimes_AA_\Fp$.

Now consider any  maximal ideal $\Fp\not=\Fp_0$ with $\Fp|a$  and set $n:=\ord_\Fp(a)$. 
By the definition of $\tTp$ we then have $\Fp^n T_\Fp(\phi) \subset \Fp^n\tTp = a\cdot M\cdot T_\Fp(\phi) \subset T_\Fp(\phi)$, and hence also $\Fp^{2n}T_\Fp(\phi) \subset \Fp^{2n}\tTp \subset \Fp^nT_\Fp(\phi)$.
But the group of $\Fp^{2n}$-division points $\phi[\Fp^{2n}](K^\sep) \cong T_\Fp(\phi)/\Fp^{2n}T_\Fp(\phi)$ and the action of $\End_{K^\sep}(\phi)$ on it can be determined by finite computation, and so can the subgroup $a\cdot M\cdot\phi[\Fp^{2n}](K^\sep) \cong \Fp^n\tTp/\Fp^{2n}T_\Fp(\phi)$. We can therefore find an explicit description of the $M_\Fp$-module $\Fp^n\tTp/\Fp^{2n}\tTp$ and its ${\End_{K^\sep}(\phi)\otimes_AA_\Fp}$-submodule $\Fp^nT_\Fp(\phi)/\Fp^{2n}\tTp$. 
Dividing by~$a$, this determines the right hand side of the cartesian diagram
\UseTheoremCounterForNextEquation
\begin{equation}\label{Cartesian}
\vcenter{\xymatrix{
\tTp \ar[r] & \tTp/\Fp^n\tTp \\
T_\Fp(\phi) \ar[r] \ar@{^{ (}->}[u] 
& T_\Fp(\phi)/\Fp^n\tTp \ar@{^{ (}->}[u] \\}}
\end{equation}
up to isomorphism. 

Note that we do not have an explicit description of~$\tTp$, but know only its isomorphism class as a projective $M_\Fp$-module. But for any positive integer $k$ the natural homomorphism 
$$\xymatrix@R-15pt{
\Aut_{M_\FP'}\bigl((M_\FP')^{\oplus k}\bigr) \ar[r] & 
\Aut_{M_\FP'}\bigl((M_\FP')^{\oplus k}/\Fp^n(M_\FP')^{\oplus k}\bigr) \\
\GL_k(M_\FP^{\prime\opp}) \ar[r] \ar@{=}[u] & 
\GL_k\bigl(M_\FP^{\prime\opp}/\Fp M_\FP^{\prime\opp}\bigr) \ar@{=}[u] \\}$$
is surjective. Thus for any finitely generated projective $M_\FP'$-module~$X$, any automorphism of $X/\Fp^nX$ lifts to an automorphism of~$X$. The same then follows for $M_\FP$ and for~$M_\Fp$; hence any automorphism of the $M_\Fp$-module $\tTp/\Fp^n\tTp$ lifts to an automorphism of the $M_\Fp$-module~$\tTp$. 
This implies that in the diagram (\ref{Cartesian}), the upper and right edges together are uniquely determined up to joint isomorphism! As the diagram is cartesian, this determines the isomorphism class of $T_\Fp(\phi)$ as a module over $\End_{K^\sep}(\phi)\otimes_AA_\Fp$, as desired.

All in all we have seen that one can effectively determine the isomorphism class of $T_\Fp(\phi)$ as a module over $\End_{K^\sep}(\phi)\otimes_AA_\Fp$ for all $\Fp\not=\Fp_0$, whence the proposition.
\end{Proof}

\begin{Thm}\label{GaloisUpToCommens}
One can effectively determine the image of the adelic Galois representation (\ref{GamadDef}) up to commensurability and conjugation under $\GL_r(A_\ad)$.
\end{Thm}

\begin{Proof}
If $\phi$ has generic characteristic, this follows by combining Theorem \ref{AllEndosKSep}, Proposition \ref{EndosOnTate}, and Theorem \ref{EndGalGen}.

If $\phi$ is non-isotrivial of special characteristic, by Propositions \ref{GoingUp} and \ref{GoingDown} one can effectively find the data $(K',\phi',f,B)$ described there, as well as an explicit presentation of $\End_{K^\sep}(\phi'|B)$. Let $\Fq_0$ be the characteristic ideal of~$\phi'|B$, and set $B_\ad := \prod_{\Fq\not=\Fq_0}B_\Fq$. Then by Proposition \ref{EndosOnTate} one can effectively determine the isomorphism class of $T_\ad(\phi'|B)$ as a module over $\End_{K^\sep}(\phi'|B)\otimes_BB_\ad$. 
% For any maximal ideal $\Fq$ of $B$ different from the characteristic ideal of $\phi'|B$, this yields an explicit description of the commutant $D_\Fq$ of $\End_{K^\sep}(\phi'|B)$ in $\End_{B_\Fq}(T_\Fq(\phi'|B))$. 
This yields an explicit description of the commutant $\prod_{\Fq\not=\Fq_0} D_\Fq$ of $\End_{K^\sep}(\phi'|B)$ in $\End_{B_\ad}(T_\ad(\phi'|B))\cong\Mat_{r''\times r''}(B_\ad)$ and hence of the group $\prod_{\Fq\not=\Fq_0} D^1_\Fq$ of elements of reduced norm~$1$.
With Theorem \ref{EndGalSpec3} one obtains a description of the image of Galois in the adelic Galois representation associated to $\phi'|B$,
up to commensurability and conjugation. The image of Galois for $\phi$ up to commensurability can be determined from this, as explained in Devic-Pink \cite[\S6.2]{DevicPink}. 

If $\phi$ is isotrivial, find $k$ and $\psi$ as in the proof of Theorem \ref{RankKSep}, so that the image of $\Gal(K^\sep/K)$ is commensurable with the pro-cyclic group generated by $\Frob_k$ associated to~$\psi$. By Proposition \ref{EndosOnTate} one can compute the action of $\Frob_k\in\End_{k^\sep}(\psi)$ on $T_\ad(\psi)\cong T_\ad(\phi)$ up to isomorphism.
% Computing the minimal polynomial of $\Frob_k$ over $A$ using Proposition \ref{MinPolEndEff} then yields the desired result.
\end{Proof}

%%%%%%%%%%%%%%%%%%%%%%%%%%%%%%%%%%%%%%%%%%%%%%%%%%%%%%%%%%%%%%%%%%%%%%%%%%%%
%\newpage
\section{Variation}
\label{Var}

In this section we briefly discuss a different approach to making the search for endomorphisms effective by bounding the degrees of generators of $\End_{K^\sep}(\phi)$ via reduction. As outlined here, this approach succeeds only in a restricted class of cases in generic characteristic, namely when $\End^{\circ}_K(\varphi)$ is separable over $F$.
% If this is known one can restrict oneself to search through the reductions of $\varphi$ and omit the search for endomorphisms. When two suitable reductions are found, one can compute the endomorphism ring directly.

\begin{Prop}\label{EndRed1}
Let $\phi$ be a Drinfeld $A$-module over an arbitrary field $L$.
Let $v$ be a valuation on~$L$ with residue field $\ell_v$ where $\phi$ has good reduction~$\phi_v$. Then the natural reduction homomorphism 
$$\End_{L}(\phi)\to\End_{\ell_v}(\phi_v)$$
is injective and the torsion of its cokernel is primary to the characteristic ideal of~$\phi_v$.
% annihilated by a power of~$\Fp$. 
\end{Prop}

\begin{Proof}
The injectivity follows from the standard fact that the degree in $\tau$ of an endomorphism is preserved under reduction.
%The injectivity is a standard result, see \cite{?}.
%Compare (\ref{EndRed})...

Let $\Fp_v$ denote the characteristic ideal of $\phi_v$. 
Extend $v$ to a valuation on $L^{\sep}$. Then the residue field of this extension is naturally a separable closure $\ell_v^{\sep}$ of $\ell_v$.
After modifying $\phi$ by an isomorphism over $L$, we can assume that it has good reduction form, meaning that $\phi$ has coefficients in the valuation ring of $v$ with highest coefficient a unit. 
Note that for any $a\in A\setminus \Fp_v$, the zeroth coefficient of $\phi_a$ is then also a unit in the valuation ring. It follows that the Newton polygon of $\phi_a(X)/X$ with respect to $v$ is a horizontal line, hence every non-zero element of $\phi[a](L^{\sep})$ has valuation zero.
	
Consider an element $f_v\in \End_{\ell_v}(\phi_v)$ and suppose there exists $a\in A\setminus \Fp_v$ such that $g_v:= f_v \phi_{v,a} $ is the reduction of some element $g\in \End_{L}(\phi)$. We claim that then $f_v$ must also lie in the image of the reduction homomorphism. 
Indeed, consider any $x\in \phi[a](L^{\sep})$ and let $x_v\in \phi_v[a](\ell_v^{\sep})$ denote its reduction. Then the reduction of $g(x)$ is $g_v(x_v)=f_v(\phi_{v,a}(x_v))=0$; hence $g(x)$ has positive valuation. Since $g$ is an endomorphism of $\phi$, we have $g(x) \in\phi[a](L^{\sep})$, and since every non-zero element of $\phi[a](L^{\sep})$ has valuation zero, it follows that $g(x)=0$. As $x$ was arbitrary, we conclude that $\phi[a](L^{\sep})\subset \Ker(g)$. Using this and the fact that $\phi_a$ is separable, we see that $g$ is right divisible by~$\phi_a$, in other words that $g=f\phi_a$ for some $f\in L[\tau]$. It is straightforward to check that this $f$ is an endomorphism of $\phi$ whose reduction is $f_v$.
	
Finally, let $h_v\in \End_{\ell_v}(\phi_v)$ be such that for some non-zero $b\in A$ the product $h_v\phi_{v,b}$ lies in the image of the reduction homomorphism. Write $bA = \Fp_v^k\Fa$ for some ideal $\Fa\subset A$ not divisible by $\Fp_v$ and pick any $a\in \Fa\setminus \Fp_v$. Then for any $c\in \Fp_v^k$, we have $b|ca$, and so $h_v\phi_{v,c}\phi_{v,a}$ lies in the image of the reduction homomorphism. Since $a\not\in \Fp_v$, by the above it follows that $h_v\phi_{v,c}$ already lies in the image. Letting $c\in \Fp_v^k$ vary, this shows that the image of $h_v$ in the cokernel of the reduction homomorphism is annihilated by $\Fp_v^k$. Letting $h_v$ vary over all elements whose image in the cokernel is torsion finishes the proof.  
\end{Proof}

\begin{Prop}\label{EndRed2}
Let $\phi$ be a Drinfeld $A$-module over an arbitrary field $L$. Let $v$ and $v'$ be valuations on~$L$ with residue fields $\ell_v$ and $\ell_{v'}$ where $\phi$ has good reduction $\phi_v$ and $\phi_{v'}$ with different characteristic ideals.
Then the image of the natural reduction homomorphism 
$$\End_{L}(\phi)\to\End_{\ell_v}(\phi_v)\times \End_{\ell_{v'}}(\phi_{v'})$$
is saturated, i.e., its cokernel is torsion free.
\end{Prop}

\begin{Proof}
Let $\Fp_v$ and $\Fp_{v'}$ be the characteristic ideals of $\phi_v$ and $\phi_{v'}$. By Proposition \ref{EndRed1}, the torsion part of the cokernel of the reduction homomorphism associated to $v$ is $\Fp_v$-primary, while the torsion part of the reduction homomorphism associated to $v'$ is $\Fp_{v'}$-primary. Using the facts that $A$ is a Dedekind ring and that the endomorphism rings are torsion free, one can show that the product $\End_L(\phi)\to \End_{\ell_v}(\phi_v)\times \End_{\ell_{v'}}(\phi_{v'}$) of the reduction homomorphisms has saturated image. This is an exercise in commutative algebra, which we leave to the reader.
	%We will show that whenever $\Fp_1$ and $\Fp_2$ are distinct maximal ideals of $A$ and for $i=1,2$ we have an injective homomorphism $j_i: M\into N_i$  of torsion free $A$-modules such that the torsion part of its cokernel is $\Fp_i$-primary, then $j_1\times j_2: M\to N_1\times N_2$ is saturated. By Proposition \ref{EndRed1}, this proves the claim. 
	%So suppose $(n_1,n_2)$ represents a torsion element in the cokernel of $j_1\times j_2$. This means there is non-zero $a\in A$ and $m\in M$ such that $a n_1 =j_1(m)$ and $a n_2 =j_2(m)$. After possibly replacing $a$ by a non-zero multiple, we can write $a= b_1 b_2$, such that $b_1\not \in\Fp_2$ and $b_2 \not\in \Fp_1$. By assumption, all zero divisors of $N_1/j_1(M)$ lie in $\Fp_1$, so from $b_2 b_1 n_1 \in j_1(M)$ obtain $b_1n_1\in j_1(M)$. Since $M,N_1$ are torsion free, it follows that $m/b_2$ must lie in $M$. Since $N_2$ is torsion free, $j_2(m/b_2)=b_1 n_2$, but all zero divisors of $N_2/j_2(M)$ lie in $\Fp_2$, so $n_2\in j_2(M)$. Using again that $N_1,N_2$ and $M$ are torsion-free, this implies $m/b_1b_2$ is in $M$ and is the preimage of $(n_1,n_2)$. This proves that $j_1\times j_2$ is saturated.	
\end{Proof}

\medskip
%We now address the assumption that $\End_K^\circ(\phi)$ is separable over $F$. First, we give some sufficient conditions that one could use in practice to decide whether this assumption is fulfilled.
Now we return to a Drinfeld $A$-module $\phi\colon A\to K[\tau]$ over a finitely generated field~$K$. As before, let $R$ be a finitely generated normal subring of $K$ such that $\phi$ extends to a Drinfeld $A$-module over $\Spec R$.

\begin{Prop}\label{OrdSep}
For any maximal ideal $\Fm$ where $\phi_\Fm$ is ordinary, $\End_{k_\Fm^\sep}^\circ(\phi_\Fm)$ is a finite separable field extension of~$F$.
\end{Prop}

\begin{Proof} 
%[REFERENCE?]
The characteristic polynomial of $\Frob_\Fm$ has precisely one root with multiplicity $1$ in $\overline{F}$ which is not a unit above the characteristic ideal of~$\phi_\Fm$. As the characteristic polynomial is a power of the minimal polynomial, it follows that the characteristic polynomial is already irreducible and separable. Thus $F(\Frob_\Fm)$ is a separable field extension of $F$ of degree $\rank(\phi_\Fm)$. 
Since $F(\Frob_\Fm)$ is the center of $\End^\circ_{k_\Fm}(\phi_\Fm)$, the formula $d_\Fm e_\Fm=\rank(\phi_\Fm)$ implies that $\End^\circ_{k_\Fm}(\phi_\Fm) = F(\Frob_\Fm)$, which is therefore commutative and separable over~$F$. The same argument over a finite extension of $k_\Fm$ proves the same for $\End_{k_\Fm^\sep}^\circ(\phi_\Fm)$.
\end{Proof}

\begin{Prop}\label{EndSep1}
If $\phi$ has generic characteristic, the following conditions are equivalent:
\begin{enumerate}
\item[(a)] $\End_K^\circ(\phi)$ is a finite separable field extension of~$F$.
\item[(b)] $\End_{K^\sep}^\circ(\phi)$ is a finite separable field extension of~$F$.
\item[(c)] There exists a maximal ideal $\Fm$ of $R$ such that the reduction $\phi_\Fm$ is ordinary.
\end{enumerate} 
Furthermore, each of the following conditions implies the ones above:
\begin{enumerate}
\item[(d)] The algebraic closure of $F$ in $K$ is separable over~$F$.
\item[(e)] The rank of $\phi$ is not divisible by~$p$.
\end{enumerate} 
\end{Prop}

\begin{Proof} 
%Since $\phi$ has generic characteristic, $\End^{\circ}_K(\phi)$ and $\End_{K^\sep}^\circ(\phi)$ are commutative and therefore finite field extensions of~$F$. As the first is contained in the second, we have (b)$\Rightarrow$(a). 
Since $\phi$ has generic characteristic, $\End^{\circ}_K(\phi)$ is commutative and thus a finite field extension of~$F$. If it is separable, the set of maximal ideals of $R$ where $\phi$ has ordinary reduction has positive Dirichlet density by \cite[Thm.\;0.3\;(b)]{PinkRIMS}. In particular it is non-empty, proving the implication (a)$\Rightarrow$(c).

Conversely, if $\phi$ has ordinary reduction at~$\Fm$, the endomorphism ring $\End^\circ_{k_\Fm}(\phi_\Fm)$ is a finite separable field extension of~$F$ by Proposition \ref{OrdSep}. The reduction of endomorphisms induces an $F$-algebra homomorphism $\End^\circ_K(\phi)\to\End^\circ_{k_\Fm}(\phi_\Fm)$. As a subfield of a separable field extension $\End_K^\circ(\phi)$ is therefore separable over~$F$. This proves the implication (c)$\Rightarrow$(a).

Thus (a) and (c) are equivalent. Since $\End_{K^\sep}^\circ(\phi) = \End_{K'}^\circ(\phi)$ for some finite separable field extension $K'$ of~$K$, and the condition (c) is invariant under extending~$K$, it follows that (a), (b), and (c) are all equivalent.

Next, in generic characteristic the natural homomorphism $\End_K(\phi)\to K$, $u=\sum_{i=0}^{d}u_i\tau^i\allowbreak \mapsto u_0$ is injective and therefore extends to an $F$-algebra homomorphism ${\End_K^{\circ}(\phi)\to K}$. Thus if the algebraic closure of $F$ in $K$ is separable over~$F$, it follows that $\End_K(\phi)$ is separable over~$F$, in other words we have (d)$\Rightarrow$(a).

Finally, in generic characteristic the degree of the field extension $\End_K^\circ(\phi)/F$ divides the rank of~$\phi$. If that is not divisible by~$p$, it follows that $\End_K^\circ(\phi)$ is separable over~$F$; in other words we have (e)$\Rightarrow$(a).
\end{Proof}
 
\begin{Lem}\label{EndSep2}
If $\phi$ has generic characteristic and $\End^\circ_K(\phi)$ is separable over $F$, then there exist maximal ideals $\Fm$ and $\Fn$ of $R$, such that the reductions $\phi_\Fm$ and $\phi_\Fn$ are ordinary and have different characteristic ideals. Furthermore, one can effectively compute ideals $\Fm$ and $\Fn$ with these properties.
\end{Lem}

\begin{Proof}
By Proposition \ref{EndSep1} there exists a maximal ideal $\Fm$ with ordinary reduction. Going through all maximal ideals of~$R$ and applying Proposition \ref{CharIsoEff}, one can therefore effectively find such an~$\Fm$.
Choose a non-zero element $s$ in the characteristic ideal $\Fp_\Fm$ of~$\phi_\Fm$. Since $\phi$ has generic characteristic, the image $s'$ of $s$ in $R$ is again non-zero. The localization $R[1/s']$ is then again a finitely generated normal subring of~$K$, and $\phi$ extends to a Drinfeld $A$-module over $\Spec R[1/s']$. Repeating the preceding argument, one can effectively find a maximal ideal of $R[1/s']$ where $\phi$ has ordinary reduction. 
Pulling this back to $R$ yields a maximal ideal $\Fn$ of $R$ where $\phi$ has ordinary reduction and whose characteristic ideal $\Fp_\Fn$ does not contain~$s$. Thus $\Fp_\Fm\not=\Fp_\Fn$, and we are done.
\end{Proof}

\medskip
If $\phi$ has generic characteristic and $\End^\circ_K(\phi)$ is separable over $F$, we can now give a different method for computing $\End_{K^\sep}(\phi)$, which does not rely on Proposition \ref{X1}. 

\begin{Thm} \label{VarAlgEff}
If $\varphi$ has generic characteristic and  $\End_{K}(\phi)$ is separable over $A$, one can effectively compute a finite separable extension $K'$ of $K$ such that $\End_{K'}(\phi)=\End_{K^\sep}(\phi)$, and a finite generating set of $\End_{K'}(\phi)$ as an $A$-module.
\end{Thm}

\begin{Proof}
%Using Proposition \ref{ModelEff} choose a normal subring $R\subset K$ which is finitely generated over $\BF_p$ with $\Quot(R)=K$ and such that $\phi$ extends to a Drinfeld $A$-module over $\Spec R$.
%Now use Proposition \ref{EndSep2} to find maximal ideals $\Fm$ and $\Fn$ of $R$ such that the reductions $\phi_\Fm$ and $\phi_\Fn$ are ordinary and have different characteristics. 
Let $\Fm$ and $\Fn$ be as in Lemma \ref{EndSep2}. Then $\End_{k_{\Fm}^{\sep}}(\phi_{\Fm})$ and $\End_{k_{\Fn}^{\sep}}(\phi_{\Fn})$ are commutative by Proposition \ref{OrdSep}.
They can be effectively computed by Proposition \ref{KFinEff}~(a) and Proposition \ref{AllEndosKFinite}. 

Compute a finite list $\CL_1$ of all $F$-subalgebras $E$ of $\End^{\circ}_{k_{\Fm}^{\sep}}(\phi_{\Fm})\times\End^{\circ}_{k_{\Fn}^{\sep}}(\phi_{\Fn})$ which are fields. This can be done by first finding all subfields of $\End^{\circ}_{k_{\Fm}^{\sep}}(\phi_{\Fm})$ containing $F$ and then determining all their $F$-embeddings into $\End^{\circ}_{k_{\Fn}^{\sep}}(\phi_{\Fn})$. By our computer algebra prerequisites both these operations can be done effectively.

For each $E$ in $\CL_1$, find finitely many generators of $E\cap \End_{k_{\Fm}^{\sep}}(\phi_{\Fm})\times\End_{k_{\Fn}^{\sep}}(\phi_{\Fn})$ as an $A$-module. Again, this is possible by our computer algebra prerequisites. Let $\CL_2$ be the set of all elements of $\End_{k_{\Fm}^{\sep}}(\phi_{\Fm})\times\End_{k_{\Fn}^{\sep}}(\phi_{\Fn})$ thus obtained as $E$ varies.

For each $x$ in $\CL_2$, compute the minimal polynomial $\min_x$ of $x$ over $F$. Choose a common splitting field $K'$ of all these over~$K$. Both steps are possible due to our computer algebra prerequisites. 
For each $x$ determine all endomorphisms of $\phi$ over $K'$ whose constant coefficient is a zero of $\min_x$, using Proposition \ref{EndodCoeff}. Let $\CL_3$ be the set of all elements of $\End_{K'}(\phi)$ thus obtained as $x$ varies.

\medskip
We claim that $\CL_3$ generates $\End_{K^\sep}(\phi)$ as an $A$-module.
To see this, note first that since $R$ is integrally closed with field of fractions $K$, the composite homomorphism $R\onto k_\Fm\into k_\Fm^\sep$ can be extended to the valuation ring associated to some valuation $v$ on~$K$. Similarly, the reduction associated to $\Fn$ comes from a valuation $v'$ on $K$. Extend $v$ and $v'$ to valuations on $K^\sep$, so that the residue fields of these extensions are $k_\Fm^\sep$ and $k_\Fn^\sep$, respectively.
%naturally separable closures $k_\Fm^\sep$ and $k_\Fn^\sep$ of $k_\Fm$ and~$k_\Fn$, respectively.
According to Proposition \ref{EndRed2}, the induced reduction homomorphism $\End_{K^\sep}(\phi)\to {\End_{k_\Fm^\sep}(\phi_\Fm)\times \End_{k_\Fn^\sep}(\phi_\Fn)}$ is injective and has saturated image.
It extends to a homomorphism $\End^\circ_{K^\sep}(\phi)\to\End^\circ_{k^\sep_\Fm}(\phi_\Fm)\times \End^\circ_{k^\sep_\Fn}(\phi_\Fn)$ of $F$-algebras, whose image $E$ is a field and must therefore appear in the list $\CL_1$.
By saturatedness, the image of $\End_{K^\sep}(\phi)$ is equal to $E \cap \End_{k^\sep_\Fm}(\phi_\Fm)\times \End_{k^\sep_\Fn}(\phi_\Fn)$. Therefore some subset $\{x_1,\ldots,x_n\}$ of $\CL_2$ will generate this image as an $A$-module. By the construction of~$\CL_3$, each $x_i$ is the reduction of an element of~$\CL_3$. This shows that $\CL_3$ generates $\End_{K^\sep}(\phi)$ as an $A$-module, as claimed.
\end{Proof}

%Other possible facts:
%
%\begin{Prop}\label{X2}
%%? ? ? (X2) For a set of places $v$ of positive density, the dimension of $\Quot(\End_{\overline{k_v}}(\phi_v))$ over its center is equal to the dimension of $\Quot(\End_{\oK}(\phi))$ over its center.
%\end{Prop}
%
%\begin{Prop}\label{X3}
%? ? ? (X3) If $Z(\End_{\oK}(\phi))=A$, there are two points $v$ and $v'$ of good reduction as in (X2), such $Z(\End_{\overline{k_v}}(\phi_v))$ and 
%$Z(\End_{\overline{k_{v'}}}(\phi_{v'}))$ are linearly disjoint over $\Quot(A)$. 
%\end{Prop}
%
%\begin{Prop}\label{X4}
%? ? ? (X4) $\Quot(\End_{\oK}(\phi))$ is the largest finite skew field extension of $F$ which possesses an embedding over $F$ into $\Quot(\End_{\overline{k_v}}(\phi_v))$ for almost all $v$.
%\end{Prop}

%%%%%%%%%%%%%%%%%%%%%%%%%%%%%%%%%%%%%%%%%%%%%%%%%%%%%%%%%%%%%%%%%%%%%%%%%%%
%\newpage
\section{Comparing two Drinfeld modules}
\label{Compare}

In this section we consider two Drinfeld $A$-modules $\phi$ and $\psi$ over $K$ with the same characteristic homomorphism $A\to K$ and hence the same characteristic ideal~$\Fp_0$. 
We will show that one can effectively decide whether $\phi$ and $\psi$ are isogenous and determine all homomorphisms, both over $\oK$ and over~$K$. 

By the Tate conjecture for $A$-motives, due to Taguchi \cite{TaguchiTate} and Tamagawa \cite{Tamagawa0} (see also Pink-Traulsen \cite[Thm.\,2.4]{PinkTraulsen2}), for any prime $\Fp\not=\Fp_0$ of $A$ we have a natural isomorphism
\UseTheoremCounterForNextEquation
\begin{equation}\label{TateHom}
\Hom_K(\phi,\psi)\otimes_AA_\Fp
\ \stackrel{\sim}{\longto}\ \Hom_{A_\Fp}(T_\Fp(\phi),T_\Fp(\psi))^{\Gal(K^\sep/K)}.
\end{equation}
In particular, as the right hand side does not change under inseparable extension, any homomorphism $\phi\to\psi$ over $\oK$ is already defined over~$K^\sep$. Thus $\phi$ and $\psi$ are isogenous over $\oK$ if and only if they are isogenous over $K^\sep$.

\begin{Prop}\label{IsogOverK}
Assume that for some prime $\Fp\not=\Fp_0$ of~$A$, all $\Fp$-torsion points of $\phi$ and $\psi$ are defined over~$K$. Then if $\phi$ and $\psi$ are isogenous over~$K^\sep$, they are isogenous over~$K$.
\end{Prop}

\begin{Proof}
Let $K'$ be a finite separable extension of $K$ over which $\phi$ and $\psi$ are isogenous. After extending $K'$ we may assume that $K'/K$ is galois. Then
$$\Hom_K(\phi,\psi)\ =\ \Hom_{K'}(\phi,\psi)^{\Gal(K'/K)}.$$
By the isomorphism (\ref{TateHom}) for $K'$ in place of~$K$ we have
$$\Hom_{K'}(\phi,\psi)\otimes_AA_\Fp
\ \stackrel{\sim}{\longto}\ \Hom_{A_\Fp}(T_\Fp(\phi),T_\Fp(\psi))^{\Gal(K^\sep/K')}.$$
The image of this isomorphism is a saturated $A_\Fp$-submodule of $\Hom_{A_\Fp}(T_\Fp(\phi),T_\Fp(\psi))$ and hence a direct summand. The induced homomorphism
$$\Hom_{K'}(\phi,\psi)\otimes_AA/\Fp
\ \longto\ \Hom_{A_\Fp}(T_\Fp(\phi),T_\Fp(\psi)) \otimes_AA/\Fp$$
is therefore injective. 
By the construction of the Tate module $T_\Fp(\phi) \otimes_AA/\Fp$ is naturally isomorphic to the group $\phi[\Fp]$ of $\Fp$-torsion points of~$\phi$, and likewise for~$\psi$. Thus we obtain a natural Galois equivariant injection
$$\Hom_{K'}(\phi,\psi)\otimes_AA/\Fp
\ \into\ \Hom_{A/\Fp}(\phi[\Fp],\psi[\Fp]).$$
By assumption $\Gal(K^\sep/K)$ acts trivially on the target group; hence $\Gal(K'/K)$ acts trivially on $\Hom_{K'}(\phi,\psi)\otimes_AA/\Fp$. Since $\Gal(K'/K)$ is a finite group, its action on $\Hom_{K'}(\phi,\psi)$ thus factors through a $p$-group and is therefore unipotent. As $\Hom_{K'}(\phi,\psi)$ is non-zero by assumption, so is consequently the submodule of $\Gal(K'/K)$-invariants. This means that $\Hom_K(\phi,\psi)$ is non-zero, as desired.
\end{Proof}

\medskip
Choose a normal integral domain $R$ that is finitely generated over~$\BF_p$ with ${\Quot(R)=K}$, such that $\phi$ and $\psi$ extend to Drinfeld $A$-modules over $\Spec R$. For any maximal ideal $\Fm$ of $R$ let $\phi_\Fm$ and $\psi_\Fm$ denote their reductions over~$k_\Fm$. 
%In place of Proposition \ref{X1} we will use the following two facts:

\begin{Prop}\label{X5}
If $\phi$ and $\psi$ are not isogenous over $K$, there exists a maximal ideal $\Fm\subset R$ such that the characteristic polynomials of $\Frob_\Fm$ associated to $\phi_\Fm$ and $\psi_\Fm$ are different.
\end{Prop}

\begin{Proof}
Suppose to the contrary that for all $\Fm$ the characteristic polynomials are equal. Pick a prime $\Fp\not=\Fp_0$ of $A$ and consider the continuous representations of $\Gal(K^\sep/K)$ on the rational Tate modules $V_\Fp(\phi) := T_\Fp(\phi)\otimes_{A_\Fp}F_\Fp$ and $V_\Fp(\psi) := T_\Fp(\psi)\otimes_{A_\Fp}F_\Fp$. As the Frobenius elements are dense in $\Gal(K^\sep/K)$, it follows that any element of $\Gal(K^\sep/K)$ has the same characteristic polynomial on $V_\Fp(\phi)$ as on $V_\Fp(\psi)$. By a general fact from representation theory (see Pink-Traulsen \cite[Prop.\,3.8]{PinkTraulsen2}) the representations thus have a common Jordan H\"older factor. But by Taguchi \cite[Thm.\,0.1]{TaguchiFinite}, \cite[Thm.\,0.1]{TaguchiInfinite} the representations are semisimple. Thus they posses an isomorphic direct summand, and in particular $\Hom_{A_\Fp}(T_\Fp(\phi),T_\Fp(\psi))^{\Gal(K^\sep/K)}$ is non-zero. By (\ref{TateHom}) it follows that $\Hom_K(\phi,\psi)$ is non-zero, contrary to the assumption.
\end{Proof}

\begin{Thm}\label{CompareKEff}
One can effectively decide whether $\phi$ and $\psi$ are isogenous over~$K$. 
\end{Thm}

\begin{Proof}
Again we start two processes in parallel:

\medskip{\it Process (a): Find isogenies:} 
For each $d\ge0$ search for isogenies $\phi\to\psi$ of degree $d$ in~$\tau$. For this choose a finite set $\CS$ of generators of the $\BF_p$-algebra~$A$. Then an element $u\in K[\tau]$ of degree~$d$ is an isogeny $\phi\to\psi$ if and only if $u\phi_a=\psi_au$ for all $a\in\CS$. With the Ansatz $u=\sum_{i=0}^{d}u_i\tau^i$ these equations amount to finitely many polynomial equations in the coefficients~$u_i$. We also know that there are at most finitely many solutions. By our computer algebra prerequisites, one can therefore effectively describe all these solutions.

As soon as an isogeny $\phi\to\psi$ is found, kill process (b) and stop with the answer ``yes''. Otherwise, repeat the calculation with $d+1$ in place of~$d$.

\medskip{\it Process (b): Compare Frobeniuses:} 
For each maximal ideal $\Fm\subset R$ use Proposition \ref{FrobEff} (a) to determine the characteristic polynomials of $\Frob_\Fm$ associated to $\phi_\Fm$ and $\psi_\Fm$. If they are different, kill process (a) and stop with the answer ``no''. Otherwise, repeat the calculation with the next~$\Fm$. 

\medskip{\it Effectivity:} 
By Proposition \ref{X5} the algorithm terminates with the correct answer.
\end{Proof}

\begin{Thm}\label{CompareKsepEff}
One can effectively decide whether $\phi$ and $\psi$ are isogenous over $K^\sep$. 
\end{Thm}

\begin{Proof}
Choose any prime $\Fp\not=\Fp_0$ of~$A$. By solving the equations for the $\Fp$-torsion points of $\phi$ and $\psi$ one can find an explicit finite separable extension $K'$ of $K$ such that all these torsion points are defined over~$K'$. Then Proposition \ref{IsogOverK} implies that $\phi$ and $\psi$ are isogenous over $K^\sep$ if and only if they are isogenous over~$K'$. This in turn can be effectively decided by Theorem \ref{CompareKEff}.
\end{Proof}

\medskip
For the remaining results we view $K^\sep[\tau]$ as an $\BF_p[t]$-module via the multiplication $(a,u)\mapsto\psi_au$.

\begin{Thm}\label{AllHomsKSep}
For any non-constant element $t\in A$, one can effectively find a finite separable extension $K''$ of $K$ and elements of $\Hom_{K''}(\phi,\psi)$ which form an orthogonal basis of $\Hom_{K^\sep}(\phi,\psi)$ over~$\BF_p[t]$.
\end{Thm}

\begin{Proof}
Use Theorem \ref{CompareKsepEff} to decide whether $\phi$ and $\psi$ are isogenous over $K^\sep$. If not, then $\Hom_{K^\sep}(\phi,\psi)=0$  with the trivial basis. If yes, the rank of $\Hom_{K^\sep}(\phi,\psi)$ over $\BF_p[t]$ is equal to that of $\End_{K^\sep}(\phi)$. The latter can be effectively determined by Theorem \ref{RankKSep} and by computing the rank of $A$ over~$\BF_p[t]$. To finish, observe that everything from Proposition \ref{OrthFindStep} through Proposition \ref{RankEff} remains true with $M_d\subset\Hom_{K^\sep}(\phi,\psi)$ in place of $\End_{K^\sep}(\phi)$ and, occasionally, $\psi$ in place of~$\phi$. The analogue of Proposition \ref{RankEff} thus yields the desired orthogonal basis. 
(But again observe Remark \ref{IsogEndEffRem}.)
\end{Proof}

\begin{Thm}\label{AllHomsK}
For any non-constant element $t\in A$, one can effectively find an orthogonal basis of $\Hom_K(\phi,\psi)$ over~$\BF_p[t]$.
\end{Thm}

\begin{Proof}
Use Theorem \ref{CompareKEff} to decide whether $\phi$ and $\psi$ are isogenous over~$K$. If not, then $\Hom_K(\phi,\psi)=0$  with the trivial basis. If yes, the rank of $\Hom_K(\phi,\psi)$ over $\BF_p[t]$ is equal to that of $\End_K(\phi)$. The latter can be effectively determined by Theorem \ref{AllEndosK}. To finish, apply the arguments from the proof of Theorem \ref{AllEndosK} to $\Hom_{K''}(\phi,\psi)$ in place of $\End_{K''}(\phi)$.
\end{Proof}

%%%%%%%%%%%%%%%%%%%%%%%%%%%%%%%%%%%%%%%%%%%%%%%%%%%%%%%%%%%%%%%%%%%%%%%%%%%

%\newpage

%%%%%%%%%%%%%%%%%%%%%%%%%%%%%%%%%%%%%%%%%%%%%%%%%%%%%%%%%%%%%%%%%%%%%%%%%%%%

\end{document}